\documentclass{article}
\usepackage{amsmath,amssymb}
\usepackage{amsthm}
\usepackage{graphicx}
\usepackage[all]{xy}
\usepackage{fancyhdr}

\theoremstyle{definition}
\newtheorem{lem}{Lemma}[section] 
\newtheorem{prop}[lem]{Proposition} 
\newtheorem{thm}[lem]{Theorem}
\newtheorem{cor}[lem]{Corollary}

\theoremstyle{definition}

\theoremstyle{remark}

\newtheorem{rem}[lem]{Remark}

\newcommand{\Dt}{D_{\tau}}
\newcommand{\Dtt}{D_{\tau}^{\times}}

\newcommand{\tih}{\tilde{h}}
\newcommand{\tiw}{\tilde{w}}

\title{A Matsumoto type theorem for $GL_n$ over rings of non-commutative Laurent polynomials}
\author{Ryusuke Sugawara}
\date{\empty}

\begin{document}
\maketitle

\footnotetext{This work was supported by JSPS KAKENHI Grant Number JP21J10690.}

\begin{abstract}
We give a Matsumoto-type presentation of $K_2$-groups over rings of
non-commutative Laurent polynomials, which is a non-commutative
version of M. Tomie's result for loop groups. Our main idea is due to
U. Rehmann's approach in case of division rings.
\end{abstract}

\section*{Introduction}

\hspace{6mm}For a field $F$, let $GL(n,F)$ and $St(n,F)$ be the general linear group of degree $n\geq 2$ over $F$ and the corresponding Steinberg group respectively. Then, $K_2(n,F)$ is defined to be the kernel of a natural homomorphism of $St(n,F)$ onto the elementary subgroup $E(n,F)$ ($=SL(n,F)$ in this case) of $GL(n,F)$.
As you can see from the Matsumoto's result \cite{hm}, such a $K_2$-group changes its group structure depending on the size of the matrix group $E(n,F)$: $K_2(2,F)$ is called of symplectic type, and presented by
the symbols $c(u,v)$, $u,v\in F^{\times}$ and the following defining relations;
\begin{align*}
&c(u,v)c(uv,w)=c(u,vw)c(v,w),\\
&c(1,1)=1,\\
&c(u,v)=c(u^{-1},v^{-1}),\\
&c(u,v)=c(u,(1-u)v)\hspace{6mm}(1-u\in F^{\times}).
\end{align*}
Also,  $K_2(n,F)$ for $n\geq 3$ is called of non-symplectic type, and presented by the symbols
$c(u,v)$, $u,v\in F^{\times}$ and the following defining relations;
\begin{align*}
&c(uv,w)=c(u,w)c(v,w),\\
&c(u,vw)=c(u,v)c(u,w),\\
&c(u,1-u)=1\hspace{6mm}(1-u\in F^{\times}).
\end{align*}
Matsumoto gave the structure of $K_2$-groups (Schur multipliers) corresponding to Chevalley groups over fields $F$, and $K_2$-groups related to Kac-Moody groups are determined by Morita and Rehmann \cite{jmur}. Especially, $K_2$-groups derived from loop groups are calculated by Tomie \cite{mt}.

Extended affine Lie algebras not included in Kac-Moody algebras are  important objects in recent study, in particular, quantum tori are remarkable ones (cf.\cite{aabgp}, \cite{yy}). Morita and Sakaguchi \cite{ms}, \cite{hs} researched groups over completed quantum tori in two variables. However, in order to deal with $n$ variables, we need rings of non-commutative Laurent polynomials over division rings different from the commutative cases $F$, $F[t,t^{-1}]$. Linear groups over division rings have already been studied by Rehmann \cite{r2} and \cite{r}, but we need more general case. We constructed a ring of non-commutative Laurent polynomials $\Dt$, and researched groups in the following exact sequence, except for a group presentation of $K_2$-groups (cf.\cite{rs}):
\[1\to K_2(n,\Dt)\to St(n,\Dt)\to GL(n,\Dt)\to K_1(n,\Dt)\to 1.\]

Let $D$ be any division ring, and fix an automorphism $\tau$ of $D$. Then $\Dt=D[t,t^{-1}]$ is a non-commutative Laurent polynomial in an indeterminate $t$, with the relation $tat^{-1}=\tau(a)$ for all $a\in D$. Let $\Dtt$ be the multiplicative group of $\Dt$, and $[\Dtt,\Dtt]$ be the commutator subgroup of $\Dtt$. In this paper, we determine the group presentation of $K_2$-groups $K_2(n,\Dt)$ over $\Dt$. Section~1 gives preliminaries: $E(n,\Dt)$ has a Bruhat decomposition controlled by the affine Weyl group $W_a$, that is, $E(n,\Dt)=\cup_{w\in W_a}BwB=UNU$. Also, we show some relations in the ``diagonal subgroup'' of  the Steinberg group. In Section~2, we determine the group presentation of symplectic $K_2$-groups $K_2(2,\Dt)$. Let $P$ be the group presented by generators $c(u,v)$, $u,v\in\Dt^{\times}$ with the following defining relations:
\begin{align*}
&\text{(P1)}\quad c(u,v)c(vu,w)=c(u,vw)c(v,w),\\
&\text{(P2)}\quad c(u,v)=c(uvu,u^{-1}),\\
&\text{(P3)}\quad c(x,y)c(u,v)c(x,y)^{-1}=c([x,y]u,v)c(v,[x,y]),\\
&\text{(P4)}\quad c(u,v)=c(u,v(1-u))\hspace{6mm}(1-u\in \Dt^{\times}),\\
&\text{(P5)}\quad c(u,v)=c(u,-vu).
\end{align*}
Then, there exists a natural homomorphism $\varphi_0$ of $P$ onto $[\Dtt,\Dtt]$, whose kernel is isomorphic to $K_2(2,\Dt)$. We will prove this along with \cite{hm}, \cite{r} and \cite{hs}. We give the presentation of non-symplectic $K_2$-groups $K_2(n,\Dt)$ in Section~3. Let $Q$ be the group presented by generators $c(u,v)$, $u,v\in\Dtt$ with the following defining relations:
\begin{align*}
&\text{(Q1)}\quad c(uv,w)=c(^uv,^uw)c(u,w),\\
&\text{(Q2)}\quad c(u,vw)=c(u,v)c(^vu,^vw),\\
&\text{(Q3)}\quad c(u,1-u)=1\hspace{6mm}(1-u\in\Dtt),
\end{align*}
where $^uv=uvu^{-1}$. Then there exists a natural homomorphism $\varphi$ of $Q$ onto $[\Dtt,\Dtt]$, whose kernel is isomorphic to $K_2(n,\Dt)$.

\begin{center}\textbf{Acknowledgments}\end{center}
\hspace{6mm}The author would like to thank Professor Jun Morita, for his continuous guidance and helpful advice.

\section{Linear groups and Steinberg groups over $\Dt$}

\hspace{5mm}Let $D$ be a division ring, and we fix an automorphism $\tau\in\mathrm{Aut}(D)$. In the following, we denote by $\Dt=D[t,t^{-1}]$ the ring of Laurent polynomials generated by $D$ and an indeterminate $t$ with the relation $tat^{-1}=\tau(a)$ for $a\in D$; note that $\Dtt=\{ st^k=t^k\tau^{-k}(s)\ |\ s\in D\setminus\{0\},k\in\mathbb{Z}\}$, where $R^{\times}$ denotes the group of invertible elements in a ring $R$. Let $M(n,\Dt)$ be the ring of $n\times n$ matrices whose entries are in $\Dt$, and define $GL(n,\Dt)=M(n,\Dt)^{\times}$.\vspace{3mm}

Let $\Delta=\{\epsilon_i-\epsilon_j\ |\ 1\leq i\neq j\leq n\}$ be the root system of type $\mathrm{A}_{n-1}$, where $\{ \epsilon_i\}_{1\leq i\leq n}$ is an orthonormal basis with respect to an inner product $(\cdot,\cdot)$. We see that any root in $\Delta$ is expressed as
\[ \epsilon_i-\epsilon_j=(\epsilon_i-\epsilon_{i+1})+(\epsilon_{i+1}-\epsilon_{i+2})+\cdots+(\epsilon_{j-1}-\epsilon_j)\]
if $i<j$ and as its minus version if $i>j$. Let $\Pi=\{ \alpha_i=\epsilon_i-\epsilon_{i+1}\ |\ 1\leq i\leq n-1\}$, be the simple system of $\Delta$. We call $\Delta^+=(\text{Span}_{\mathbb{Z}_{\geq 0}}\Pi)\cap\Delta$ a set of positive roots in $\Delta$, and $\Delta^-=-\Delta^+$ a set of negative roots in $\Delta$. Also, let $\Delta_a=\Delta\times\mathbb{Z}$ be the (abstract) affine root system of type $\mathrm{A}_{n-1}^{(1)}$, and let $\Pi_a=\{ \dot{\alpha_i}=(\alpha_i,0)\ |\ 1\leq i\leq n-1\}\cup\{\dot{\alpha_0}=(-\theta,1)\}$ be the simple system of $\Delta_a$, where $\theta=\alpha_1+\dots+\alpha_{n-1}$ is a highest root in $\Delta$. We denote by $\Delta_a^+=(\Delta^+\times\mathbb{Z}_{\geq 0})\cup(\Delta^-\times\mathbb{Z}_{>0})$ (resp. $\Delta_a^-=(\Delta^+\times\mathbb{Z}_{< 0})\cup(\Delta^-\times\mathbb{Z}_{\leq 0})$) be the set of positive roots (resp. negative roots) in $\Delta_a$. We identify $\beta\in \Delta$ with $(\beta,0)\in \Delta_a$, and regard $\Delta$ as a subset of $\Delta_a$. For $\dot{\beta}=(\beta,m)\in\Delta_a$, we set $-\dot{\beta}=(-\beta,-m)$.\vspace{3mm}

In a standard way,  the Weyl group $W$ of $\Delta$ is generated by all reflections $\sigma_{\beta}$ for $\beta\in\Delta$, and the Weyl group $W_a$ of $\Delta_a$ (the affine Weyl group of $\Delta$) is generated by all reflections $\sigma_{\dot{\beta}}$ for $\dot{\beta}\in\Delta_a$. The action of $\sigma_{\dot{\beta}}$ is defined by
\[\sigma_{\dot{\beta}}(\dot{\gamma})=(\sigma_{\beta}(\gamma), n-\langle \gamma, \beta \rangle m)\]
for $\dot{\beta}=(\beta, m),\ \dot{\gamma}=(\gamma, n)\in\Delta_a$ and $\langle \gamma,\beta\rangle=2(\gamma,\beta)/(\beta,\beta)$.\vspace{3mm}

For $\beta=\epsilon_i-\epsilon_j\in\Delta$, $f\in \Dt$ we define
	\[ x_{\beta}(f)=x_{ij}(f)=I+fE_{ij}\]
where $I\in GL(n,\Dt)$ is the identity matrix and $E_{ij}\in M(n,\Dt)$ is the matrix unit. For $\beta=\epsilon_i-\epsilon_j\in\Delta$ and $u\in\Dtt$, we set
\begin{align*}
	w_{\beta}(u)&=w_{ij}(u)=x_{\beta}(u)x_{-\beta}(-u^{-1})x_{\beta}(u),\\
	h_{\beta}(u)&=h_{ij}(u)=w_{\beta}(u)w_{\beta}(-1).
\end{align*}
Also, for $\dot{\beta}=(\beta, m)\in\Delta_a$, $f\in D$, and $s\in D^{\times}$, we set
\begin{align*}
	& x_{\dot{\beta}}(f)=\begin{cases}
						x_{\beta}(ft^m) & \text{if $\beta\in\Delta^+$}, \\
						x_{\beta}(t^mf) & \text{if $\beta\in\Delta^-$},
					 \end{cases}\\
	&w_{\dot{\beta}}(s)=x_{\dot{\beta}}(s)x_{-\dot{\beta}}(-s^{-1})x_{\dot{\beta}}(s),\\
	&h_{\dot{\beta}}(s)=w_{\dot{\beta}}(s)w_{\beta}(-1).
\end{align*}
It can be easily checked that $x_{\dot{\beta}}(f)^{-1}=x_{\dot{\beta}}(-f)$, $w_{\dot{\beta}}(u)^{-1}=w_{\dot{\beta}}(-u)$. The elementary subgroup $E(n,\Dt)$ is defined to be the subgroup of $GL(n,\Dt)$ generated by $x_{\dot{\beta}}(f)$ for all $\dot{\beta}\in\Delta_a$ and $f\in D$, that is,
\begin{align*}
E(n,\Dt)
	&=\langle x_{\dot{\beta}}(f)\ |\ \dot{\beta}=(\beta,m)\in\Delta_a,\ f\in D\rangle\\
	&=\langle x_{ij}(g)\ |\ 1\leq i\neq j\leq n,\ g\in\Dt\rangle.
\end{align*}

\noindent As subgroups of $E(n,\Dt)$, we put
\begin{align*}
	&U_{\dot{\beta}}=\{ x_{\dot{\beta}}(f)\ |\ f\in D\}\ \text{for $\dot{\beta}\in \Delta_a$},\\
	&U'_{\pm\dot{\alpha}}=\langle x_{\pm\dot{\alpha}}(g)U_{\dot{\beta}}x_{\pm\dot{\alpha}}(g)^{-1}\ |\ g\in D,\ \dot{\beta}\in\Delta^{\pm}_a\setminus\{ \pm\dot{\alpha}\}\rangle\quad\text{for $\dot{\alpha}\in\Pi_a$},\\
	&U^{\pm}=\langle U_{\dot{\beta}}\ |\ \dot{\beta}\in\Delta_a^{\pm}\rangle,\\
	&N=\langle w_{\dot{\beta}}(u)\ |\ \dot{\beta}\in\Delta_a,\ u\in D^{\times}\rangle,\\
	&T=\langle h_{\dot{\beta}}(u)\ |\ \dot{\beta}\in\Delta_a,\ u\in D^{\times}\rangle.
\end{align*}
Notice that an element $h\in T$ is of the form $h=\mathrm{diag}(u_1,u_2,\dots,u_n)$ with some $u_i\in\Dtt$. For each $1\leq i\leq n$, if $u_i=s_it^{m_i}$ with $s_i\in D^{\times}$ and $m_i\in\mathbb{Z}$, then we define
	\[ \mathrm{deg}_i(h)=\mathrm{deg}(u_i)=m_i. \]
We set
\begin{align*}
	&T_0=\langle h\ |\ h\in T,\ \mathrm{deg}_i(h)=0\ \text{for all}\ i=1,2,\dots,n\rangle,\\
	&B^{\pm}=\langle U^{\pm},\ T_0\rangle,\\
	&S=\{ w_{\dot{\beta}}(1)\ \mathrm{mod}\ T_0\ |\ \dot{\beta}\in\Delta_a\}.
\end{align*}

\noindent We can easily check two ``elementary'' relations in $E(n,\Dt)$ as follows (see \cite{rs}): for $\beta=\epsilon_i-\epsilon_j$, $\gamma=\epsilon_k-\epsilon_l\in\Delta$ with $i\neq l,j\neq k$, and $f,g\in D_{\tau}$,
\begin{align*}
&\text{(R1)}\ x_{\beta}(f)x_{\beta}(g)=x_{\beta}(f+g),\\
&\text{(R2)}\ [x_{\beta}(f),x_{\gamma}(g)]=
	\begin{cases}
		x_{\beta+\gamma}(fg) \hspace{5mm}&\text{if }\ \beta+\gamma\in\Delta,\ j=k,\\
		x_{\beta+\gamma}(-gf) \hspace{5mm}&\text{if }\ \beta+\gamma\in\Delta,\ i=l,\\
		I \hspace{10mm} &\text{otherwise}.
	\end{cases}
\end{align*}
We also give relations between $x_{\beta}(f)$, $w_{\beta}(u)$, and $h_{\beta}(u)$. For $\beta=\epsilon_i-\epsilon_j, \gamma=\epsilon_k-\epsilon_l\in\Delta$, $f,g\in D_{\tau}$, and $s,u\in D_{\tau}^{\times}$,
\begin{align*}
&\text{(R3)}\ w_{\beta}(u)x_{\gamma}(f)w_{\beta}(u)^{-1}=
	\begin{cases}
		x_{\gamma}(f)\hspace{5mm}&\text{if} \hspace{2mm}(\beta,\gamma)=0,\\
		x_{\mp\beta}(-u^{\mp 1}fu^{\mp 1})\hspace{5mm}&\text{if} \hspace{2mm}\gamma=\pm\beta,\\
		x_{\sigma_{\beta}(\gamma)}(-u^{-1}f)\hspace{5mm}&\text{if} \hspace{2mm}\beta\pm\gamma\neq 0\ \text{and}\ i=k,\\
		x_{\sigma_{\beta}(\gamma)}(-fu)\hspace{5mm}&\text{if} \hspace{2mm}\beta\pm\gamma\neq 0\ \text{and}\ i=l,\\
		x_{\sigma_{\beta}(\gamma)}(uf)\hspace{5mm}&\text{if} \hspace{2mm}\beta\pm\gamma\neq 0\ \text{and}\ j=k,\\
		x_{\sigma_{\beta}(\gamma)}(fu^{-1})\hspace{5mm}&\text{if} \hspace{2mm}\beta\pm\gamma\neq 0\ \text{and}\ j=l,\\
	\end{cases}\\
&\text{(R4)}\ h_{\beta}(u)x_{\gamma}(f)h_{\beta}(u)^{-1}=
	\begin{cases}
		x_{\gamma}(f)\hspace{5mm}&\text{if} \hspace{2mm}(\beta,\gamma)=0,\\
		x_{\pm\beta}(-u^{\pm 1}fu^{\pm 1})\hspace{5mm}&\text{if} \hspace{2mm}\gamma=\pm\beta,\\
		x_{\gamma}(uf)\hspace{5mm}&\text{if} \hspace{2mm}\beta\pm\gamma\neq 0\ \text{and}\ i=k,\\
		x_{\gamma}(fu^{-1})\hspace{5mm}&\text{if} \hspace{2mm}\beta\pm\gamma\neq 0\ \text{and}\ i=l,\\
		x_{\gamma}(u^{-1}f)\hspace{5mm}&\text{if} \hspace{2mm}\beta\pm\gamma\neq 0\ \text{and}\ j=k,\\
		x_{\gamma}(fu)\hspace{5mm}&\text{if} \hspace{2mm}\beta\pm\gamma\neq 0\ \text{and}\ j=l,\\
	\end{cases}\\
&\text{(R5)}\ w_{\beta}(u)w_{\gamma}(s)w_{\beta}(u)^{-1}=
	\begin{cases}
		w_{\gamma}(s)\hspace{5mm}&\text{if} \hspace{2mm}(\beta,\gamma)=0,\\
		w_{\mp\beta}(-u^{\mp 1}su^{\mp 1})\hspace{5mm}&\text{if} \hspace{2mm}\gamma=\pm\beta,\\
		w_{\sigma_{\beta}(\gamma)}(-u^{-1}s)\hspace{5mm}&\text{if} \hspace{2mm}\beta\pm\gamma\neq 0\ \text{and}\ i=k,\\
		w_{\sigma_{\beta}(\gamma)}(-su)\hspace{5mm}&\text{if} \hspace{2mm}\beta\pm\gamma\neq 0\ \text{and}\ i=l,\\
		w_{\sigma_{\beta}(\gamma)}(us)\hspace{5mm}&\text{if} \hspace{2mm}\beta\pm\gamma\neq 0\ \text{and}\ j=k,\\
		w_{\sigma_{\beta}(\gamma)}(su^{-1})\hspace{5mm}&\text{if} \hspace{2mm}\beta\pm\gamma\neq 0\ \text{and}\ j=l,\\
	\end{cases}\\
&\text{(R6)}\ w_{\beta}(u)h_{\gamma}(s)w_{\beta}(u)^{-1}=
	\begin{cases}
		h_{\gamma}(s)\hspace{5mm}&\text{if} \hspace{2mm}(\beta,\gamma)=0,\\
		h_{\mp\beta}(u^{\mp 1}su^{\mp 1})h_{\mp\beta}(u^{\pm 2})\hspace{5mm}&\text{if} \hspace{2mm}\gamma=\pm\beta,\\
		h_{\sigma_{\beta}(\gamma)}(u^{-1}s)h_{\sigma_{\beta}(\gamma)}(u)\hspace{5mm}&\text{if} \hspace{2mm}\beta\pm\gamma\neq 0\ \text{and}\ i=k,\\
		h_{\sigma_{\beta}(\gamma)}(su)h_{\sigma_{\beta}(\gamma)}(u^{-1})\hspace{5mm}&\text{if} \hspace{2mm}\beta\pm\gamma\neq 0\ \text{and}\ i=l,\\
		h_{\sigma_{\beta}(\gamma)}(us)h_{\sigma_{\beta}(\gamma)}(u^{-1})\hspace{5mm}&\text{if} \hspace{2mm}\beta\pm\gamma\neq 0\ \text{and}\ j=k,\\
		h_{\sigma_{\beta}(\gamma)}(su^{-1})h_{\sigma_{\beta}(\gamma)}(u)\hspace{5mm}&\text{if} \hspace{2mm}\beta\pm\gamma\neq 0\ \text{and}\ j=l.\\
	\end{cases}
\end{align*}

We know from \cite[Theorem 1]{rs} the following proposition.\vspace{3mm}

{\prop\label{prop:11}{ $(E(n,\Dt),B,N,S)$ is a Tits system with the corresponding affine Weyl group $W_a$. In particular, $E(n,\Dt)$ has a Bruhat decomposition $E(n,\Dt)=\cup_{w\in W_a}BwB=UNU$.}}\vspace{3mm}

Let $St(n,\Dt)$ be the Steinberg group, which is defined by the generators $\hat{x}_{\beta}(f)=\hat{x}_{ij}(f)$ for $f\in\Dt$, $\beta=\epsilon_i-\epsilon_j\in\Delta$, and the following defining relations:
\begin{align*}
&\text{(ST1)}\ \hat{x}_{\beta}(f)\hat{x}_{\beta}(g)=\hat{x}_{\beta}(f+g),\\
&\text{(ST2)}\ [\hat{x}_{\beta}(f),\hat{x}_{\gamma}(g)]=
		\begin{cases}
			\hat{x}_{\beta+\gamma}(fg) &\text{if $\beta+\gamma\in\Delta$ and $j=k$},\\
			\hat{x}_{\beta+\gamma}(-gf) &\text{if $\beta+\gamma\in\Delta$ and $i=l$}, \\
			1 &\text{otherwise},
		\end{cases}
\end{align*}
where $f,g\in\Dt$ and $\beta=\epsilon_i-\epsilon_j,\gamma=\epsilon_k-\epsilon_l\in\Delta$ with $\gamma\neq \pm\beta$. Then, we can easily check that $\hat{x}_{\beta}(f)^{-1}=\hat{x}_{\beta}(-f)$ from (ST1). When $n=2$, we use
\[\text{(ST2)}'\ \hat{w}_{\beta}(u)\hat{x}_{\beta}(f)\hat{w}_{\beta}(-u)=\hat{x}_{-\beta}(-u^{-1}fu^{-1})\]
instead of (ST2), where $\hat{w}_{\beta}(u)=\hat{x}_{\beta}(u)\hat{x}_{-\beta}(-u^{-1})\hat{x}_{\beta}(u)$ for $u\in\Dtt$. Then there exists a natural homomorphism $\phi$ of $St(n,\Dt)$ onto $E(n,\Dt)$
defined by $\phi(\hat{x}_{ij}(f))=x_{ij}(f)$ for $1\leq i\neq j\leq n$ and $f\in \Dt$.\vspace{3mm}

For $\dot{\beta}=(\beta,m)\in\Delta_a$ and $f\in D$, we set
\[\hat{x}_{\dot{\beta}}(f)=
				 \begin{cases}
					\hat{x}_{\beta}(ft^m) &\text{if $\beta\in\Delta^+$}, \\
					\hat{x}_{\beta}(t^mf) &\text{if $\beta\in\Delta^-$},
				 \end{cases}\]
and we also put $\hat{w}_{\dot{\beta}}(u)=\hat{x}_{\dot{\beta}}(u)\hat{x}_{-\dot{\beta}}(-u^{-1})\hat{x}_{\dot{\beta}}(u)$
and $\hat{h}_{\dot{\beta}}(u)=\hat{w}_{\dot{\beta}}(u)\hat{w}_{\beta}(-1)$ for $u\in D^{\times}$. We put, as subgroups of $St(n,\Dt)$,
\begin{align*}
	&\hat{U}_{\dot{\beta}}=\{ \hat{x}_{\dot{\beta}}(f)\ |\ f\in D\},\\
	&\hat{U}^{\pm}=\langle \hat{U}_{\dot{\beta}}\ |\ \dot{\beta}\in\Delta_a^{\pm}\rangle,\\
	&\hat{N}=\langle \hat{w}_{\dot{\beta}}(u)\ |\ \dot{\beta}\in\Delta_a,\ u\in D^{\times}\rangle,\\
	&\hat{T}=\langle \hat{h}_{\dot{\beta}}(u)\ |\ \dot{\beta}\in\Delta_a,\ u\in D^{\times}\rangle.
\end{align*}
Also, we set
\begin{align*}
	&\hat{T_0}=\langle h\ |\ h\in \hat{T},\ \mathrm{deg}_i(\phi(h))=0\ \mathrm{for\ all}\ i=1,\dots, n\rangle,\\
	&\hat{B}^{\pm}=\langle \hat{U}^{\pm},\ \hat{T_0}\rangle,\\
	&\hat{S}=\{ \hat{w}_{\dot{\beta}}(1)\ \mathrm{mod}\ \hat{T_0}\ |\ \dot{\beta}\in\Delta_a\}.
\end{align*}

We give several relations between $\hat{x}_{\dot{\beta}}(f)$, $\hat{w}_{\dot{\beta}}(u)$ and $\hat{h}_{\dot{\beta}}(s)$ like (R1)--(R6) in $St(n,\Dt)$, but those are the same except (R6). In the Steinberg group $St(n,\Dt)$, we get ($\hat{R}$6) as follows: for $\beta=\epsilon_i-\epsilon_j$, $\gamma=\epsilon_k-\epsilon_l\in\Delta$, and $u,s\in\Dt^{\times}$,

\begin{align*}
(\hat{R}6)\ &\hat{w}_{\beta}(u)\hat{h}_{\gamma}(s)\hat{w}_{\beta}(u)^{-1}\\&=
	\begin{cases}
		\hat{h}_{\gamma}(s)\hspace{2mm}&\text{if} \hspace{2mm}(\beta,\gamma)=0,\\
		\hat{h}_{\mp\beta}(-u^{\mp 1}su^{\mp 1})\hat{h}_{\mp\beta}(-u^{\pm 2})^{-1} &\text{if} \hspace{2mm}\gamma=\pm\beta,\\
		\hat{h}_{\sigma_{\beta}(\gamma)}(-u^{-1}s)\hat{h}_{\sigma_{\beta}(\gamma)}(-u^{-1})^{-1} &\text{if} \hspace{2mm}\beta\pm\gamma\neq 0\ \text{and}\ i=k,\\
		\hat{h}_{\sigma_{\beta}(\gamma)}(-su)\hat{h}_{\sigma_{\beta}(\gamma)}(-u)^{-1} &\text{if} \hspace{2mm}\beta\pm\gamma\neq 0\ \text{and}\ i=l,\\
		\hat{h}_{\sigma_{\beta}(\gamma)}(us)\hat{h}_{\sigma_{\beta}(\gamma)}(u)^{-1} &\text{if} \hspace{2mm}\beta\pm\gamma\neq 0\ \text{and}\ j=k,\\
		\hat{h}_{\sigma_{\beta}(\gamma)}(su^{-1})\hat{h}_{\sigma_{\beta}(\gamma)}(u^{-1})^{-1} &\text{if} \hspace{2mm}\beta\pm\gamma\neq 0\ \text{and}\ j=l.
	\end{cases}
\end{align*}

{\prop[cf.\cite{rs}]\label{prop:12}{$(St(n,\Dt),\hat{B}^{\pm},\hat{N},\hat{S})$ is a Tits system with the corresponding affine Weyl group $W_a$.}}\vspace{3mm}

Then, we know from \cite[Theorems 3 and 5]{rs} that $\phi$ is a central extension, in particular,
a universal central extension if $|Z(D)|\geq 5$ and $|Z(D)|\neq 9$, where $Z(D)$ is the center of $D$. Here, we define our $K_2$-groups as $K_2(n,\Dt)=\mathrm{Ker}\phi$ (cf.\cite{m}). Then we get the following proposition from \cite[Proposition 5]{rs}:\vspace{3mm}

{\prop\label{prop:13}{Let $n\geq 2$. We have
\begin{align*}
K_2(n,\Dt)=\langle \hat{c}(u_1,v_1)\hat{c}(u_2,v_2)&\dotsm\hat{c}(u_r,v_r) | u_i,v_i\in\Dt^{\times},\\
&[u_1,v_1][u_2,v_2]\dotsm [u_r,v_r]=1\rangle,
\end{align*}
where $\hat{c}(u,v)=\hat{h}_{12}(u)\hat{h}_{12}(v)\hat{h}_{12}(vu)^{-1}$, and
$K_2(n,\Dt)$ is a central subgroup of $St(n,\Dt)$.}}\vspace{3mm}

We give some relations in $\hat{T}$ for later use (see Section 2). We first assume that $n=2$. Then we get the next lemma (cf.\cite{r}).\vspace{3mm}

{\lem\label{lem:14}{In $\hat{T}$, the followings hold for $1\leq i\neq j\leq 2$ and $u,v,x,y\in D_{\tau}^{\times}$.}}
\begin{align*}
&\mathrm{(T1)}\quad \hat{h}_{ij}(u)\hat{h}_{ij}(v)=\hat{h}_{ij}(uvu)\hat{h}_{ij}(u^{-1}),\\
&\mathrm{(T1)'}\quad \hat{h}_{ij}(u)\hat{h}_{ij}(v)=\hat{h}_{ij}(v^{-1})\hat{h}_{ij}(vuv),\\
&\mathrm{(T2)}\quad \hat{c}(u,v)=\hat{h}_{ij}(v^{-1}u^{-1})^{-1}\hat{h}_{ij}(u^{-1})\hat{h}_{ij}(v^{-1}),\\
&\mathrm{(T2)'}\quad \hat{h}_{ij}(u^{-1}v^{-1})^{-1}\hat{h}_{ij}(u^{-1})\hat{h}_{ij}(v^{-1})=\hat{h}_{ij}(u)\hat{h}_{ij}(v)\hat{h}_{ij}(uv)^{-1},\\
&\mathrm{(T3)}\quad \hat{c}_{ij}(u,v)=\hat{c}_{ij}(uvu,u^{-1})=\hat{c}_{ij}(v^{-1},vuv),\\
&\mathrm{(T4)}\quad \hat{c}_{ij}(u,v)\hat{c}_{ij}(vu,w)=\hat{c}_{ij}(u,vw)\hat{c}_{ij}(v,w)\\
	&\hspace{35.3mm}=\hat{h}_{ij}(u)\hat{c}_{ij}(v,w)\hat{h}_{ij}(u)^{-1}\hat{c}_{ij}(u,wv),\\
&\mathrm{(T5)}\quad \hat{c}_{ij}(x,y)\hat{h}_{ij}(u)\hat{c}_{ij}(x,y)^{-1}=\hat{h}_{ij}([x,y]u)\hat{h}_{ij}([x,y])^{-1}\\
	&\hspace{46mm}=\hat{h}_{ij}([y,x])^{-1}\hat{h}_{ij}(u[y,x]),\\
&\mathrm{(T6)}\quad \hat{c}_{ij}(x,y)\hat{c}_{ij}(u,v)\hat{c}_{ij}(x,y)^{-1}=\hat{c}_{ij}(u,[x,y])^{-1}\hat{c}_{ij}(u,[x,y]v)\\
	&\hspace{48.5mm}=\hat{c}_{ij}([x,y]u,v)\hat{c}_{ij}([x,y],v)^{-1},\\
&\mathrm{(T7)}\quad \hat{c}_{ij}(u,v)=\hat{c}_{ij}(u,v(1-u))\hspace{5mm}\text{if $1-u\in\Dt^{\times}$},\\
&\mathrm{(T8)}\quad \hat{c}_{ij}(u,v)=\hat{c}_{ij}(u,-vu).
\end{align*}

\begin{proof}
(T1) and (T1)$'$ can be shown by direct calculation. Let us show (T2) and (T2)$'$.
By (T1), we see that
\begin{align*}
&\hat{h}_{ij}(u)\hat{h}_{ij}(u^{-1}v)=\hat{h}_{ij}(vu)\hat{h}_{ij}(u^{-1}),\\
&\hat{h}_{ij}(v)\hat{h}_{ij}(v^{-1}u^{-1})=\hat{h}_{ij}(u^{-1}v)\hat{h}_{ij}(v^{-1}),
\end{align*}
and hence $\hat{h}_{ij}(vu)^{-1}\hat{h}_{ij}(v)\hat{h}_{ij}(u)=\hat{h}_{ij}(u^{-1})\hat{h}_{ij}(v^{-1})\hat{h}_{ij}(v^{-1}u^{-1})^{-1}$.
Replacing $u$ and $v$ by $u^{-1}$ and $v^{-1}$, respectively, we obtain
\[\hat{h}_{ij}(v^{-1}u^{-1})^{-1}\hat{h}_{ij}(v^{-1})\hat{h}_{ij}(u^{-1})=\hat{h}_{ij}(u)\hat{h}_{ij}(v)\hat{h}_{ij}(vu)^{-1}=\hat{c}_{ij}(u,v),\]
as desired. (T3) is a straightforward consequence of (T1) and (T1)$'$.
We show the first equality in (T4) as follows:
\begin{align*}
&\hat{c}_{ij}(u,v)\hat{c}_{ij}(vu,w)\\
&=\hat{h}_{ij}(u)\hat{h}_{ij}(v)\hat{h}_{ij}(w)\hat{h}_{ij}(wvu)^{-1}\text{(by the definitions of $\hat{c}_{ij}(u,v)$ and $\hat{c}_{ij}(vu,w)$)}\\
&=\hat{c}_{ij}(u,vw)\hat{h}_{ij}(vwu)\hat{h}_{ij}(vw)^{-1}\hat{h}_{ij}(v)\hat{h}_{ij}(w)\hat{h}_{ij}(wvu)^{-1}\\
&=\hat{c}_{ij}(u,vw)\hat{h}_{ij}(vwu)\hat{h}_{ij}(v^{-1})\hat{h}_{ij}(w^{-1})\hat{h}_{ij}(v^{-1}w^{-1})^{-1}\hat{h}_{ij}(wvu)^{-1}\\
&=\hat{c}_{ij}(u,vw)\hat{h}_{ij}(v)\hat{h}_{ij}(wuv^{-1})\hat{h}_{ij}(v^{-1})^{-1}\hat{h}_{ij}(v^{-1})&\text{(by (T1)$'$)}\\
&\hspace{30mm}\times\hat{h}_{ij}(w^{-1})\hat{h}_{ij}(v^{-1}w^{-1})^{-1}\hat{h}_{ij}(wvu)^{-1}\\
&=\hat{c}_{ij}(u,vw)\hat{h}_{ij}(v)\hat{h}_{ij}(w)\hat{h}_{ij}(uv^{-1}w^{-1})\hat{h}_{ij}(w^{-1})^{-1}&\text{(by (T1))}\\
&\hspace{30mm}\times\hat{h}_{ij}(w^{-1})\hat{h}_{ij}(v^{-1}w^{-1})^{-1}\hat{h}_{ij}(wvu)^{-1}\\
&=\hat{c}_{ij}(u,vw)\hat{c}_{ij}(v,w)\hat{h}_{ij}(wv)\hat{h}_{ij}(uv^{-1}w^{-1})\hat{h}_{ij}(v^{-1}w^{-1})^{-1}\hat{h}_{ij}(wvu)^{-1}\\
&=\hat{c}_{ij}(u,vw)\hat{c}_{ij}(v,w).&\text{(by (T1))}
\end{align*}
For the second equality in (T4), we compute
\begin{align*}
\hat{c}_{ij}(u,v)\hat{c}_{ij}(vu,w)&=\hat{h}_{ij}(u)\hat{h}_{ij}(v)\hat{h}_{ij}(w)\hat{h}_{ij}(wvu)^{-1}\\
	&=\hat{h}_{ij}(u)\hat{c}_{ij}(v,w)\hat{h}_{ij}(u)^{-1}\hat{h}_{ij}(u)\hat{h}_{ij}(wv)\hat{h}_{ij}(wvu)^{-1}\\
	&=\hat{h}_{ij}(u)\hat{c}_{ij}(v,w)\hat{h}_{ij}(u)^{-1}\hat{c}_{ij}(u,wv).
\end{align*}
For (T5), we need some extra identities. Putting $v=1$ in (T1) and (T1)$'$, we find that
\[ \hat{h}_{ij}(u)=\hat{h}_{ij}(u^2)\hat{h}_{ij}(u^{-1})=\hat{h}_{ij}(u^{-1})\hat{h}_{ij}(u^2),\]
which implies that $\hat{h}_{ij}(u^2)^{-1}=\hat{h}_{ij}(u^{-2})$. Thus we get
\[\text{(T5)}'\ \hat{h}_{ij}(u)^{\pm 1}\hat{h}_{ij}(v)\hat{h}_{ij}(u)^{\mp 1}=\hat{h}_{ij}(u^{\pm 1}vu^{\pm 1})\hat{h}_{ij}(u^{\pm 2})^{-1}.\]
Applying (T5)$'$ three times, we obtain the first equality in (T5) as follows:
\begin{align*}
&\hat{c}_{ij}(x,y)\hat{h}_{ij}(u)\hat{c}_{ij}(x,y)^{-1}\\
&=\hat{h}_{ij}(x)\hat{h}_{ij}(y)\hat{h}_{ij}(yx)^{-1}\hat{h}_{ij}(u)\hat{h}_{ij}(yx)\hat{h}_{ij}(v)^{-1}\hat{h}_{ij}(x)^{-1}\\
&=\hat{h}_{ij}(x)\hat{h}_{ij}(y)\hat{h}_{ij}(x^{-1}y^{-1}ux^{-1}y^{-1})\hat{h}_{ij}(x^{-1}y^{-1}x^{-1}y^{-1})^{-1}\hat{h}_{ij}(y)^{-1}\hat{h}_{ij}(x)^{-1}\\
&=\hat{h}_{ij}(x)\hat{h}_{ij}(yx^{-1}y^{-1}ux^{-1})\hat{h}_{ij}(yx^{-1}y^{-1}x^{-1})^{-1}\hat{h}_{ij}(x)^{-1}\\
&=\hat{h}_{ij}([x,y]u)\hat{h}_{ij}([x,y])^{-1}.
\end{align*}
Also, by using (T1), we get the second equality in (T5) as follows:
\begin{align*}
&\hat{c}_{ij}(x,y)\hat{h}_{ij}(u)\hat{c}_{ij}(x,y)^{-1}\\
&=\hat{h}_{ij}([x,y]u)\hat{h}_{ij}([x,y])^{-1}\\
&=\hat{h}_{ij}([y,x])^{-1}\hat{h}_{ij}(u[y,x])\hat{h}_{ij}(u[y,x])^{-1}\hat{h}_{ij}([y,x])\hat{h}_{ij}([x,y]u)\hat{h}_{ij}([x,y])^{-1}\\
&=\hat{h}_{ij}([y,x])^{-1}\hat{h}_{ij}(u[y,x])\hat{h}_{ij}(u[y,x])^{-1}\hat{h}_{ij}(u[y,x])\\
&=\hat{h}_{ij}([y,x])^{-1}\hat{h}_{ij}(u[y,x]).
\end{align*}
(T6) is a consequence of (T4) and (T5). (T7) is equivalent to
\[\hat{w}_{ij}(vu-vu^2)\hat{w}_{ij}(v-vu)^{-1}=\hat{w}_{ij}(vu)\hat{w}_{ij}(v)^{-1},\]
which follows from
\begin{align*}
\hat{w}_{ij}&(vu-vu^2)\hat{w}_{ij}(v-vu)^{-1}\\
&=\hat{x}_{ij}(vu^2)\hat{x}_{ij}(-vu^2)\hat{w}_{ij}(vu-vu^2)\hat{w}_{ij}(v-vu)^{-1}\\
&=\hat{x}_{ij}(vu^2)\hat{w}_{ij}(vu(1-u))\hat{w}_{ij}(-v(1-u))\hat{x}_{ij}(-v)&\text{(by (ST2))}\\
&=\hat{x}_{ij}(vu)\hat{x}_{ji}(-(1-u)^{-1}u^{-1}v^{-1})\hat{w}_{ji}((1-u)^{-1}v^{-1})\\
&\hspace{30mm}\times\hat{x}_{ji}(-(1-u)^{-1}uv^{-1})\hat{x}_{ij}(-v)&\text{(by (R1))}\\
&=\hat{x}_{ij}(vu)\hat{x}_{ji}(-u^{-1}v^{-1})\hat{x}_{ij}(vu-u)\hat{x}_{ji}(v^{-1})\hat{x}_{ij}(-v)&\text{(by (R1))}\\
&=\hat{w}_{ij}(vu)\hat{w}_{ij}(v)^{-1}.
\end{align*}
(T8) is equivalent to
\[\hat{w}_{ij}(vu)\hat{w}_{ij}(v)\hat{w}_{ij}(-vu)=\hat{w}_{ij}(vu^2),\]
and which follows from (R3) and $\hat{w}_{ij}(u)=\hat{w}_{ji}(-u^{-1})$. Thus we have proved Lemma \ref{lem:14}.
\end{proof}\vspace{3mm}

Next, assume that $n\geq 3$. We know from \cite[Lemma 9.8, Lemma 9.10]{m} that 
\begin{align*}
&\text{(TT1)}\quad \hat{h}_{ij}(u)\hat{h}_{ji}(u)=1,\\
&\text{(TT2)}\quad \hat{h}_{ij}(u)\hat{h}_{ki}(u)\hat{h}_{jk}(u)=1&\text{($k\neq i$, $k\neq j$)},\\
&\text{(TT3)}\quad \hat{c}_{ij}(s,1-s)=1
\end{align*}
for $1\leq i\neq j \leq n$, $u\in\Dtt$, and $s\in D^{\times}$ with $s\neq 1$. Also, we know from \cite[Proposition 2.1]{r2} that
\[\text{(TT0)}\quad \text{$\hat{c}_{ij}(u,v)=\hat{c}_{ik}(u,v)$ for $k\neq i$ and $u,v\in\Dtt$}.\]

{\lem\label{lem:15}{ The following relations hold in $\hat{T}$.
\begin{align*}
&\text{(TT4)}\quad \hat{c}_{ij}(u,v)=[\hat{h}_{ij}(u),\hat{h}_{ik}(v)],\\
&\text{(TT4)$'$}\quad \hat{c}_{ij}(u,v)^{-1}=\hat{c}_{ij}(v,u),\\
&\text{(TT5)}\quad \hat{h}_{ij}(x)\hat{c}_{ik}(u,v)\hat{h}_{ij}(u)^{-1}\\
	&\hspace{15mm}=\hat{c}_{ik}(u,x)^{-1}\hat{c}_{ik}(u,xv)
	=\hat{c}_{ik}(xu,v)\hat{c}_{ik}(x,v)^{-1}
	=\hat{c}_{ik}(^xu,^xv),\\
&\text{(TT6)}\quad \hat{c}_{ij}(uv,w)=\hat{c}_{ij}(^uv,^uw)\hat{c}_{ij}(u,w),\\
&\text{(TT7)}\quad \hat{c}_{ij}(u,vw)=\hat{c}_{ij}(u,v)\hat{c}_{ij}(^vu,^vw).
\end{align*}}}
\begin{proof}
(TT4) can be easily checked by ($\hat{R}6$). (TT4)$'$ is computed as follows by (TT0) and (TT4):
\[\hat{c}_{ik}(u,v)^{-1}=[\hat{h}_{ik}(u),\hat{h}_{ij}(v)]^{-1}=\hat{c}_{ij}(v,u)=\hat{c}_{ik}(v,u).\]
By ($\hat{R}$6), we see that
\begin{align*}
\hat{h}_{ij}(x)\hat{c}_{ik}(u,v)\hat{h}_{ij}(u)^{-1}
&=\hat{h}_{ij}(x)\hat{h}_{ik}(u)\hat{h}_{ik}(v)\hat{h}_{ik}(vu)^{-1}\hat{h}_{ij}(x)^{-1}\\
&=\hat{h}_{ik}(xu)\hat{h}_{ik}(x)^{-1}\hat{h}_{ik}(xv)\hat{h}_{ik}(xvu)^{-1}\\
&=\hat{c}_{ik}(u,x)^{-1}\hat{c}_{ik}(u,xv).
\end{align*}
The first equation of (TT5) holds. Similarly, we get
\[\hat{h}_{ij}(x)\hat{c}_{ik}(u,v)\hat{h}_{ij}(u)^{-1}=\hat{c}_{ik}(xu,v)\hat{c}_{ik}(x,v)^{-1}\]
by (TT4)$'$. In order to prove the last equality of (TT5), we need the next relation; for $u,v,w\in\Dtt$,
\[\text{(TT5)}'\quad \hat{c}_{ij}(u,vw)=\hat{c}_{ij}(uv,w)\hat{c}_{ij}(wu,v).\]
We compute
\begin{align*}
&\hat{c}_{ij}(u,vw)\\
&=[\hat{h}_{ij}(u),\hat{h}_{ik}(vw)]\\
&=[\hat{h}_{ij}(u),\hat{c}_{ik}(v,w)\hat{h}_{ik}(w)\hat{h}_{ik}(v)]\\
&=\hat{h}_{ij}(u)\hat{c}_{ik}(v,w)\hat{h}_{ik}(w)\hat{h}_{ik}(v)\hat{h}_{ij}(u)^{-1}\hat{h}_{ik}(v)^{-1}\hat{h}_{ik}(w)^{-1}\hat{c}_{ik}(v,w)^{-1}\\
&=\hat{h}_{ij}(u)\hat{c}_{ik}(v,w)\hat{h}_{ik}(w)\hat{h}_{ik}(v)\hat{h}_{ij}(u)^{-1}\hat{h}_{ik}(vw)^{-1}\\
&=\hat{c}_{ik}(uv,w)\hat{c}_{ik}(u,w)^{-1}\hat{h}_{ij}(u)\hat{h}_{ik}(w)\hat{h}_{ik}(v)\hat{h}_{ij}(u)^{-1}\hat{h}_{ik}(vw)^{-1}\quad\text{(by (TT5))}\\
&=\hat{c}_{ik}(uv,w)\hat{h}_{ik}(w)\hat{h}_{ij}(u)\hat{h}_{ik}(v)\hat{h}_{ij}(u)^{-1}\hat{h}_{ik}(vw)^{-1}\\
&=\hat{c}_{ik}(uv,w)\hat{h}_{ik}(w)\hat{c}_{ij}(u,v)\hat{h}_{ik}(v)\hat{h}_{ik}(vw)^{-1}\\
&=\hat{c}_{ik}(uv,w)\hat{c}_{ij}(wu,v)\hat{c}_{ij}(w,v)^{-1}\hat{h}_{ik}(w)\hat{h}_{ik}(v)\hat{h}_{ik}(vw)^{-1}\quad\text{(by (TT5))}\\
&=\hat{c}_{ij}(uv,w)\hat{c}_{ij}(wu,v)\quad\text{(by (TT0))}.
\end{align*}
Therefore, we obtain
\begin{align*}
\hat{c}_{ik}(^xu,^xv)&=\hat{c}_{ik}(xu,vx^{-1})\hat{c}_{ik}(vux^{-1},x)\quad\text{(by (TT5)$'$)}\\
	&=\hat{c}_{ik}(xu,v)\hat{h}_{ij}(v)\hat{c}_{ik}(xu,x^{-1})\hat{c}_{ik}(ux^{-1},x)\hat{h}_{ij}(v)\hat{c}_{ik}(v,x)\quad\text{(by (TT5))}\\
	&=\hat{c}_{ik}(xu,v)\hat{c}_{ik}(v,x)\quad\text{(by (TT5)$'$)}\\
	&=\hat{h}_{ij}(x)\hat{c}_{ik}(u,v)\hat{h}_{ij}(x)^{-1}\quad\text{(by (TT5))}.
\end{align*}
Hence (TT5) holds. (TT6) and (TT7) are shown by (TT5).
\end{proof}\vspace{3mm}


\section{Symplectic $K_2$-groups}

\hspace{5mm}Let $D_{\tau}=D[t,t^{-1}]$ be the ring of non-commutative Laurent polynomials defined in Section~1. Let $P$ be the group generated by $c(u,v)$, $u,v\in D_{\tau}^{\times}$, with the defining relations that for $u,v,w,x,y\in D_{\tau}^{\times}$ and $s\in D^{\times}$ with $s\neq 1$,
\begin{align*}
&\text{(P1)}\ c(u,v)c(vu,w)=c(u,vw)c(v,w),\\
&\text{(P2)}\ c(u,v)=c(uvu,u^{-1}),\\
&\text{(P3)}\ c(x,y)c(u,v)c(x,y)^{-1}=c([x,y]u,v)c(v,[x,y]),\\
&\text{(P4)}\ c(s,v)=c(s,v(1-s)),\\
&\text{(P5)}\ c(u,v)=c(u,-vu).
\end{align*}

\noindent We can show the following proposition by applying (P5) to (P1).
{\lem\label{lem:21}{ In $P$, it holds that for $u,v\in D_{\tau}$,}
\begin{align*}
\text{(P5)}'\ c(u,v)=c(-uv,v).
\end{align*}}

\noindent Because $[u,v]\in [D_{\tau}^{\times},D_{\tau}^{\times}]$, $u,v\in D_{\tau}^{\times}$, satisfy the same relations (P1)--(P5) (with $c(u,v)$ replaced by $[u,v]$), there exists a (unique) surjective group homomorphism $\varphi_0: P\twoheadrightarrow [D_{\tau}^{\times},D_{\tau}^{\times}]$ which sends $c(u,v)$ to $[u,v]$ for $u,v\in D_{\tau}^{\times}$. Set $L_0=\mathrm{Ker}\ \varphi_0$; note that
\begin{align*}
L_0=\{ c(u_1,v_1)^{p_1}c(u_2,v_2)^{p_2}\cdots &c(u_r,v_r)^{p_r}\ |\ p_i=\pm 1,\ u_i,v_i\in\Dt^{\times},\\
&[u_1,v_1]^{p_1}[u_2,v_2]^{p_2}\cdots [u_r,v_r]^{p_r}=1\}.&\text{(\#)}
\end{align*}
We can show the following lemma in the same manner as \cite[Lemma 2.2]{r}.\vspace{3mm}

{\prop\label{prop:22}{ The following exact sequence is a central extension of $P$ by $L_0$.
\[1\longrightarrow L_0\longrightarrow P\stackrel{\varphi_0}{\longrightarrow} [D_{\tau}^{\times},D_{\tau}^{\times}]\longrightarrow 1.\]
}}\vspace{3mm}

By comparing (P1), (P2), (P3), (P4), (P5) with (T4), (T3), (T6), (T7), (T8), recpectively, we see that there exists a (unique) surjective group homomorphism $\zeta_0: P\to K_2(2,D_{\tau})$ with maps $c(u,v)$ to $\hat{c}(u,v)$ for $u,v\in D_{\tau}^{\times}$. By Proposition \ref{prop:13} and ($\#$), the restriction of this $\zeta_0$ to $L_0\subset P$ is a surjective group homomorphism from $L_0$ onto $K_2(2,D_{\tau})$. The following is the main theorem of this paper.\vspace{3mm}

{\thm\label{thm:23}{
The group homomorphism $\zeta_0: L_0 \rightarrow K_{2}(2,D_{\tau})$ is
an isomorphism of groups.
}}\vspace{3mm}

We proceed the proof of this theorem along with \cite{hm}, \cite{r} and \cite{hs}. \vspace{3mm}

Let $\tilde{H}_0$ be the group generated by $\tilde{h}(u)$ for $u\in\Dt^{\times}$ and $z(l)$ for $l\in P$ with the defining relations that
\begin{align*}
&\text{(H1)}\ \tilde{h}(u)\tilde{h}(v)=\tilde{h}(uvu)\tilde{h}(u^{-1}),\\
&\text{(H2)}\ \tilde{h}(u)\tilde{h}(v)=z(c(u,v))\tilde{h}(vu),\\
&\text{(H3)}\ z(l_1)z(l_2)=z(l_1l_2),\\
&\text{(H4)}\ \tilde{h}(u)z(l)=z(c(u,\varphi_0(l)))z(l)\tilde{h}(u)
\end{align*} 
for $u,v\in\Dt^{\times}$ and $l,l_1,l_2\in P$. We deduce by (H3) that $\{z(l)\ |\ l\in P\}$ is a subgroup of $\tilde{H}_0$. Moreover, we see by a way similar to \cite[CHAPITRE II]{hm} and \cite[Proposition 2]{jmur} that $L_0$ is isomophic to $\{z(l)\ |\ l\in L_0\}$. We identify $\{z(l)\ |\ l\in L_0\}$ with $L_0$, and write $z(l)$ simply by $l$ for $l\in L_0$.\vspace{3mm}

We can easily check the next proposition. 
{\prop\label{prop:24}{ All relations in Lemma \ref{lem:14} with $\hat{h}_{12}(u)$ replaced by $\tilde{h}(u)$ and $\hat{c}(u,v)$ replaced by $z(c(u,v))$ for $u,v\in D_{\tau}^{\times}$ hold in $\tilde{H}_0$.
}}\vspace{3mm}

{\prop\label{prop:25}{There exists a (unique) surjective group homomorphism $\pi_0: \tilde{H}_0\to T$ which sends $\tilde{h}(u)$ to $h_{12}(u)$ for $u\in D_{\tau}^{\times}$, and $c(v,w)$ to $h_{12}(v)h_{12}(w)h_{12}(wv)^{-1}$ for $v,w\in D_{\tau}^{\times}$. The kernel $\mathrm{Ker}\ \pi_0$ is identical to $L_0$. Moreover, 
\[1\longrightarrow L_0\longrightarrow \tilde{H}_0\stackrel{\pi_0}{\longrightarrow} T\longrightarrow 1\]
 is a central extension of $T$ by $L_0$.}}\vspace{3mm}
\begin{proof}
The first assertion is obvious by $h_{21}(u)=h_{12}(u^{-1})\in T$ and Proposition \ref{prop:24}. We first deduce, in exactly the same way as \cite[STEP 1]{rs}, that every element $h$ in $\tilde{H}_0$ can be written in the form $h=\xi\tilde{h}(s)$ for some $s\in D_{\tau}^{\times}$, where $\xi=z(c(u_1,v_1)^{p_1}c(u_2,v_2)^{p_2}\cdots c(u_r,v_r)^{p_r})\in \{z(l)\ |\ l\in P\}$ with $u_i,v_i\in D_{\tau}^{\times}$ and $p_i\in\{\pm 1\}$. Therefore, if $h\in\mathrm{Ker}\ \pi_0$, then we have $s=1$, hence $\tilde{h}(s)=\tilde{h}(1)=1$ by (T5)$'$. By induction on $r$, we see that
\[ h\tilde{h}(u)h^{-1}=\xi\tilde{h}(u)\xi^{-1}=\tilde{h}(\pi_0(\xi)u)\tilde{h}(\pi_0(\xi))^{-1}.\]
By $\pi_0(\xi)=1$, $h$ is central in $\tilde{H}_0$.
\end{proof}\vspace{3mm}

\noindent By Propositions \ref{prop:22} and \ref{prop:25}, we get the commutative diagram;
\begin{equation}
\begin{split}
\xymatrix{
&L_0 \ar[r]\ar@{=}[d] &P \ar[r]^-{\varphi_0}\ar[d]^-{z}  &[D_{\tau}^{\times},D_{\tau}^{\times}]\ar[d]^-{d}\\
&L_0 \ar[r] &\tilde{H}_0 \ar[r]^-{\pi_0}  &T,
}
\end{split}
\tag{CD}
\end{equation}
where $d$ is the embedding defined by $d([u,v])=\mathrm{diag}([u,v],1)$. This implies that $P$ is isomorphic to $\{z(l)\ |\ l\in P\}$. We identify $\{z(l)\ |\ l\in P\}$ with $P$, and write $z(l)$ simply by $l$ for $l\in P$.

{\lem\label{lem:26}{
If $\tilde{h}\in\tilde{H}_0$ is such that $\pi_0(\tilde{h})=\left(\begin{array}{cc}u_1&0\\ 0&u_2\end{array}\right)$ with some $u_1,u_2\in D_{\tau}^{\times}$, then for $u\in D_{\tau}^{\times}$, 
\[\tilde{h}\tilde{h}(u)\tilde{h}^{-1}=\tilde{h}(u_1uu_2^{-1})\tilde{h}(u_1u_2^{-1})^{-1}.\]}}
\begin{proof}
As mentioned in the proof of Proposition \ref{prop:25}, an element $\tilde{h}\in \tilde{H}_0$ can be written in the form $\tilde{h}=\tilde{h}(v)\xi$ with some $v\in\Dt^{\times}$ and $\xi\in P$ by (H4). If we put $\varphi_0(\xi)=s\in [\Dt^{\times},\Dt^{\times}]$, then $\pi_0(\tilde{h})=\mathrm{diag}(vs,v^{-1})$, and
\begin{align*}
\tilde{h}\tilde{h}(u)\tilde{h}^{-1}&=\tilde{h}(v)\tilde{h}(su)\tilde{h}(s)^{-1}\tilde{h}(v)^{-1}\\
&=\tilde{h}(vsuv)\tilde{h}(vsv)^{-1}
\end{align*}
by (T1) and (T5).
\end{proof}\vspace{3mm}

Next, we construct some extension of the monomial subgroup $N$, which is ``compatible'' with the extension $(\tilde{H}_0,\pi_0)$ of $T$ in Proposition \ref{prop:25} (see Proposition \ref{prop:29} below). For this, we give the presentation of $N$, and then define an action of $N$ on $\tilde{H}_0$. The next lemma and proposition follow from \cite[Proposition 3]{jmur} and \cite[Proposition 5.8]{hs}.\vspace{3mm}

{\lem\label{lem:27}{
The subgroup $N$ of $E(2,D_{\tau})$ is the group generated by $w_{12}(u)$ for $u\in \Dt^{\times}$ with the defining relations that for $u,v\in D_{\tau}^{\times}$,}
\begin{align*}
&\text{(N1)}\ w_{12}(u)^{-1}=w_{12}(-u),\\
&\text{(N2)}\ w_{12}(1)h_{12}(u)w_{12}(1)^{-1}=h_{12}(u^{-1}),\\
&\text{(N3)}\ h_{12}(u)h_{12}(v)=h_{12}(uvu)h_{12}(u^{-1}).
\end{align*}

{\prop\label{prop:28}{
There exists an action of $N$ on $\tilde{H}_0$ defined by
\[w_{12}(u)\cdot \tilde{h}(v)=\tilde{h}(uv^{-1}u)\tilde{h}(u^2)^{-1}\]
for $u,v\in\Dt^{\times}$; we denote $w_{12}(u)\cdot \tilde{h}(v)$ also by $w_{12}(u)\tilde{h}(v)w_{12}(u)^{-1}$ for convenience.}}\vspace{3mm}

Let $\tilde{W}_0=\langle\vartheta\rangle\cong \mathbb{Z}$ be the cyclic group of infinite order. We put $\tilde{T}_0=\langle\vartheta^2 \rangle\subset \tilde{W}_0$, and $N^*_0=\tilde{W}_0\ltimes\tilde{H}_0$, where $\tilde{W}_0$ acts on $\tilde{H}_0$ by
\[ \vartheta\cdot\tilde{h}(u)=w_{12}(-1)\cdot\tilde{h}(u)\quad\text{for $u\in D_{\tau}^{\times}$}.\]
There exists a group homomorphism $\eta_0:\tilde{T}_0\to\tilde{H}_0$ defined by $\eta_0(\vartheta^2)=\tilde{h}(-1)$. Let $J_0$ be the normal subgroup of $N^*_0=\tilde{W}_0\ltimes\tilde{H}_0$ generated by $(\vartheta^{2n},\eta_0(\vartheta^{2n})^{-1})$ for $n\in\mathbb{Z}$, and let $\tilde{N}_0=N^*_0/J_0$ be the quotient group. Let $w_1=\vartheta J_0$ be the coset in $\tilde{N}_0$ containing $\vartheta\in\tilde{W}_0\subset N^*_0$. Let $\tilde{\psi}_0:N^*_0\twoheadrightarrow N^*_0/J_0=\tilde{N}_0$ be the canonical homomorphism of $N^*_0$. We deduce that the restriction of $\tilde{\psi}_0$ to $\tilde{H}_0\subset N^*_0$ is injective; we regard $\tilde{H}_0$ as a subgroup of $\tilde{N}_0$. If we set $\tilde{w}(u)=\tilde{h}(u)\vartheta^{-1}J_0\in\tilde{N}_0$ for $u\in D_{\tau}^{\times}$, then in $\tilde{N}_0=N^*_0/J_0$,
\begin{align*}
\tilde{w}(-1)&=\tilde{h}(-1)\vartheta^{-1}J_0\\
&=\tilde{h}(-1)\tilde{h}(-1)^{-1}\vartheta J_0&\text{(by $\tilde{h}(-1)=\vartheta^2$)}\\
&=\vartheta J_0,\vspace{3mm}\\
\tilde{w}(u)\tilde{w}(-u)&=\tilde{h}(u)\vartheta^{-1}\tilde{h}(-u)\vartheta^{-1}J_0\\
	&=\tilde{h}(u)(\vartheta^{-1}\cdot\tilde{h}(-u))\vartheta^{-2}J_0\\
	&=\tilde{h}(u)\tilde{h}(-u^{-1})\tilde{h}(-1)^{-1}J_0&\text{(by Proposition \ref{prop:28})}\\
	&=c(u,-u^{-1})J_0&\text{(by the definition of $c(\cdot,\cdot)$)}\\
	&=c(1,u^{-1})J_0&(\text{by (P5)})\\
	&=J_0.
\end{align*}
Therefore, it holds that $\tilde{w}(-1)=\vartheta$, $\tilde{w}(u)^{-1}=\tilde{w}(-u)$, and $\tilde{h}(u)=\tilde{w}(u)\tilde{w}(-1)$ in $\tilde{N}_0$. Notice that there exists a group homomorphism $\psi^*_0$ from $N^*_0=\tilde{W}_0\ltimes\tilde{H}_0$ to $N$ such that $\psi^*_0(\vartheta)=w_{12}(-1)$ and $\psi^*_0(\tilde{h})=\pi_0(\tilde{h})$ for $\tilde{h}\in\tilde{H}_0$ since it can easily checked that the relations in $N^*_0$ hold in $N$. Since $\psi^*_0(\vartheta^2 \eta_0(\vartheta^2)^{-1})=w_{12}(-1)^2h_{12}(-1)^{-1}=1$, we see that $J_0\subset\mathrm{Ker}\ \psi^*_0$. Let $\psi_0:\tilde{N}_0\to N$ be the induced group homomorphism which sends $\tilde{w}(u)$ to $w_{12}(u)$. Then we see that $\psi_0$ is surjective by $w_{21}(u)=w_{12}(-u^{-1})\in N$, and is a central extension of $N$ by $L_0$. Using this and Proposition \ref{prop:25}, we deduce that the restriction of $\psi_0$ to $\tilde{H}_0$ is injective (see (CD)). We obtain, in consequence, the following proposition.\vspace{3mm}

{\prop\label{prop:29}{ The kernel $\mathrm{Ker}\ \psi_0$ of the group homomorphism $\psi_0: \tilde{N}_0\to N$ is contained in the center of $\tilde{N}_0$, and is isomorphic to $L_0$. Namely,
\[1\longrightarrow L_0\longrightarrow \tilde{N}_0\stackrel{\psi_0}{\longrightarrow} N\longrightarrow 1\]
is a central extension of $N$ by $L_0$. Moreover, the restriction of $\psi_0$ to $\tilde{H}_0$ coincides with the group homomorphism $\pi_0:\tilde{H}_0\to T\subset N$ defined in Proposition \ref{prop:25}.}}

{\lem\label{lem:210}{
If $h\in\tilde{H}_0\subset\tilde{N}_0$ is such that $\psi_0(h)=\left(\begin{array}{cc}u_1&0\\ 0&u_2\end{array}\right)\in T$ with some $u_1,u_2\in D_{\tau}^{\times}$, then
\[hw_1h^{-1}=\tilde{h}(u_1u_2^{-1})w_1.\]}}

\begin{proof}
We first show that
\begin{align*}
\text{(L1)}\ w_1c(u,v)w_1^{-1}&=\tilde{h}([u,v])^{-1}c(u,v),\\
\text{(L2)}\ w_1c(u,v)w_1^{-1}&=\tilde{h}([v,u])c(u,v).
\end{align*}
We see by (T1) and (T1)$'$ that
\begin{align*}
&\tilde{h}([u,v])w_1c(u,v)w_1^{-1}\\
&=\tilde{h}(uvu^{-1}v^{-1})\tilde{h}(u^{-1})\tilde{h}(v^{-1})\tilde{h}(u^{-1}v^{-1})^{-1}&\text{(by Proposition \ref{prop:28})}\\
&=\tilde{h}(u)\tilde{h}(v)\tilde{h}(u^{-1}v^{-1}u^{-1}v^{-1})\tilde{h}(u^{-1}v^{-1})^{-1}&\text{(by (T1)$'$)}\\
&=c(u,v)\tilde{h}(vu)\tilde{h}(u^{-1}v^{-1}u^{-1}v^{-1})\tilde{h}(u^{-1}v^{-1})^{-1}\\
&=c(u,v).&\text{(by (T1))}
\end{align*}
In a way similar to (L1), we get (L2). These equalities imply that
\begin{align*}
&\text{(L3)}\ \tilde{h}(s)^{-1}=\tilde{h}(s^{-1})\quad\text{for $s\in [\Dt^{\times},\Dt^{\times}]$},\\
&\text{(L4)}\ w_1\xi w_1^{-1}=\tilde{h}(\psi_0(\xi)^{-1})\xi\quad\text{for $\xi\in P$}.
\end{align*}
As in the proof of Lemma \ref{lem:26}, we write $h\in\tilde{H}_0$ as $h=\tilde{h}(v)\xi$ with some $v\in D_{\tau}^{\times}$ and $\xi\in P$. If we set $\psi_0(\xi)=s\in [\Dt^{\times},\Dt^{\times}]$, then $\psi_0(h)=\begin{pmatrix}vs&0\\ 0&v^{-1}\end{pmatrix}$. Thus we compute
\begin{align*}
\tilde{h}w_1\tilde{h}^{-1}&=\tilde{h}(v)\xi w_1\xi^{-1}\tilde{h}(v)^{-1}\\
&=\tilde{h}(v)\tilde{h}(s^{-1})^{-1}w_1\tilde{h}(v)^{-1}&\text{(by (L4))}\\
&=\tilde{h}(v)\tilde{h}(s^{-1})^{-1}(w_1\cdot\tilde{h}(v)^{-1})w_1\\
&=\tilde{h}(v)\tilde{h}(s)\tilde{h}(u^{-1})^{-1}w_1&\text{(by (L3) and Proposition \ref{prop:28})}\\
&=\tilde{h}(vsv)w_1.&\text{(by (T1))}
\end{align*}
\end{proof}

We know the following lemma from \cite[Lemma 5.10]{hs}.

{\lem\label{lem:211}{ Every matrix $e\in E(2,\Dt)$ can be written as $e=uwv$ with some $u,v\in U$ and $w\in N$. Moreover, the monomial matrix part $w$ is uniquely determined by $e$; we define $\rho_0: E(2,\Dt)\to N$ by $\rho_0(e)=\rho(uwv)=w$.}\vspace{3mm}

Next, we determine the value of $\rho_0$ after multiplying a double coset by $w_{\dot{\alpha_1}}(\pm 1)$ and $w_{\dot{\alpha_0}}(\pm 1)$. We remark that each $w\in N$ is either of the following forms:
\begin{equation}
\underbrace{w=\left(\begin{array}{cc}
u_1&0\\
0&u_2
\end{array}\right)}_{\text{diagonal}}\quad\text{or}\quad
\underbrace{w=\left(\begin{array}{cc}
0&u_2\\
u_1&0
\end{array}\right)}_{\text{anti-diagonal}}\tag{$\spadesuit$}
\end{equation}
with some $u_1,u_2\in D_{\tau}^{\times}$. Here we set $d=d_w=\mathrm{deg}(u_1^{-1}u_2)$. Then, for $\dot{\beta}=(\beta,m)\in\Delta_a$, $\beta=\epsilon_i-\epsilon_j$ and $f\in D$, we get the following equations.\vspace{3mm}

\underline{Case 1} If $w$ is diagonal, then
\begin{align*}
w^{\pm 1}x_{\dot{\beta}}(f)w^{\mp 1}=
	\begin{cases}
		x_{(\beta,m\mp d)}(\tau^{m\mp d}(t^{\pm d}\tau^{-m}(u_1^{\pm 1}f)u_2^{\mp 1}))&\text{if }\beta\in\Delta^+,\\
		x_{(\beta,m\mp d)}(\tau^{-m\pm d}(u_2^{\pm 1}\tau^m(fu_1^{\mp 1})t^{\pm d}))&\text{if }\beta\in\Delta^-.
	\end{cases}
\end{align*}

\underline{Case 2} If $w$ is anti-diagonal, then
\begin{align*}
w^{\pm 1}x_{\dot{\beta}}(f)w^{\mp 1}=
	\begin{cases}
		x_{(-\beta,m-d)}(\tau^{-m+d}(u_{\pm}^{\pm 1}f\tau^m(u_{\mp}^{\mp 1})t^d)&\text{if }\beta\in\Delta^+,\\
		x_{(-\beta,m-d)}(\tau^{m-d}(t^d\tau^{-m}(u_{\mp}^{\pm 1})fu_{\pm}^{\mp 1}))&\text{if }\beta\in\Delta^-,
	\end{cases}\\
\text{where $u_+=u_1$ and $u_-=u_2$.}
\end{align*}

In what follows, if $wx_{\dot{\beta}}(f)w^{-1}=x_{\dot{\gamma}}(g)$ for suitable elements $w\in N$, $f,g\in D$, and $\dot{\beta},\dot{\gamma}\in\Delta_a$, then we denote $\dot{\gamma}$ by $w(\dot{\beta})$. Here, we know from \cite[Proposition 2]{rs} that every element $e\in E(2,D_{\tau})$ can be written in the form $e=yx_{{\dot{a}}}(-f)wx_{{\dot{b}}}(g)z$ with some $f,g\in D$, ${\dot{a}},{\dot{b}}\in\Pi_a$, $y\in U'_{\dot{a}}$, and $z\in U'_{\dot{b}}$.\vspace{3mm}

{\lem\label{lem:212}{
For $e\in E(2,\Dt)$ let $\rho_0(e)=w$ be as in Lemma \ref{lem:211}, and set $e=yx_{{\dot{a}}}(-f)wx_{{\dot{b}}}(g)z$ for $f,g\in D$, ${\dot{a}},{\dot{b}}\in\Pi_a$, $y,\in U'_{\dot{a}}$ and $z\in U'_{\dot{b}}$.
Then we obtain:\vspace{3mm}

\underline{Case 1} (for $w_{\dot{a}}(1)e$).\vspace{2mm}

If $f=0$ or $w^{-1}(\dot{a})\in\Delta_a^+$, then $\rho_0(w_{\dot{a}}(1)e)=w_{\dot{a}}(1)w$.

If $f\neq 0$ and $w^{-1}(\dot{a})\not\in\Delta_a^+$, then $\rho_0(w_{\dot{a}}(1)e)=h_{\dot{a}}(1)h_{\dot{a}}(f)^{-1}w$.\vspace{3mm}

\underline{Case 2} (for $ew_{\dot{b}}(-1)$).\vspace{2mm}

If $g=0$ or $w(\dot{b})\in\Delta_a^+$, then $\rho_0(ew_{\dot{b}}(-1))=ww_{\dot{b}}(-1)$.

If $g\neq 0$ and $w(\dot{b})\not\in\Delta_a^+$, then $\rho_0(ew_{\dot{b}}(-1))=wh_{\dot{b}}(g)h_{\dot{b}}(1)^{-1}$.}}\vspace{3mm}

\begin{proof}
It suffices to calculate that
\[w_{\dot{a}}(1)x_{\dot{a}}(-f)w\hspace{5mm}\mathrm{and}\hspace{5mm}wx_{\dot{b}}(g)w_{\dot{b}}(-1).\]

If $f=0$ or $w^{-1}(\dot{a})\in\Delta_a^+$, then
\[w_{\dot{a}}(1)x_{\dot{a}}(-f)w=\underbrace{w_{\dot{a}}(1)w}_{N}\underbrace{x_{w^{-1}(\dot{a})}(k_1)}_{U},\]
where
\begin{align*}
&k_1=\begin{cases}
	-u_1^{-1}fu_2t^{-d}&\text{if $\dot{a}=\dot{\alpha_1}$ and $w$ is diagonal},\\
	-t^du_2^{-1}fu_1&\text{if $\dot{a}=\dot{\alpha_1}$ and $w$ is not diagonal},\\
	-t^{-d}\tau^{-1}(u_2^{-1})fu_1&\text{if $\dot{a}=\dot{\alpha_0}$ and $w$ is diagonal},\\
	-u_1^{-1}\tau(fu_2)t^d&\text{if $\dot{a}=\dot{\alpha_0}$ and $w$ is not diagonal}.
\end{cases}
\end{align*}

If $f\neq 0$ and $w^{-1}(\dot{a})\not\in\Delta_a^+$, then
\[w_{\dot{a}}(1)x_{\dot{a}}(-f)w=\underbrace{x_{\dot{a}}(f^{-1})}_{U}\underbrace{h_{\dot{a}}(1)h_{\dot{a}}(f)^{-1}w}_{N}\underbrace{x_{w^{-1}(-\dot{a})}(k_2)}_{U},\]
where
\begin{align*}
&k_2=\begin{cases}
	-t^{-d}u_2^{-1}f^{-1}u_1&\text{if $\dot{a}=\dot{\alpha_1}$ and $w$ is diagonal},\\
	-u_1^{-1}f^{-1}u_2t^d&\text{if $\dot{a}=\dot{\alpha_1}$ and $w$ is not diagonal},\\
	-u_1^{-1}f^{-1}\tau^{-1}(u_2)t^{-d}&\text{if $\dot{a}=\dot{\alpha_0}$ and $w$ is diagonal},\\
	-t^d\tau(u_2^{-1}f^{-1})u_1&\text{if $\dot{a}=\dot{\alpha_0}$ and $w$ is not diagonal}.
\end{cases}
\end{align*}

If $g=0$ or $w(\dot{b})\in\Delta_a^+$, then
\[wx_{\dot{b}}(g)w_{\dot{b}}(-1)=\underbrace{x_{w(\dot{b})}(k_3)}_{U}\underbrace{ww_{\dot{b}}(-1)}_{N},\]
where 
\begin{align*}
&k_3=\begin{cases}
	u_1gu_2^{-1}t^d&\text{if $\dot{b}=\dot{\alpha_1}$ and $w$ is diagonal},\\
	t^du_1gu_2^{-1}&\text{if $\dot{b}=\dot{\alpha_1}$ and $w$ is not diagonal},\\
	t^d\tau^{-1}(u_2)gu_1^{-1}&\text{if $\dot{b}=\dot{\alpha_0}$ and $w$ is diagonal},\\
	u_2\tau(gu_1^{-1})t^d&\text{if $\dot{b}=\dot{\alpha_0}$ and $w$ is not diagonal}.
\end{cases}
\end{align*}

If $g\neq 0$ and $w(\dot{b})\not\in\Delta_a^+$, then
\[wx_{\dot{b}}(g)w_{\dot{b}}(-1)=\underbrace{x_{w(-\dot{b})}(k_4)}_{U}\underbrace{wh_{\dot{b}}(g)h_{\dot{b}}(1)^{-1}}_{N}\underbrace{x_{\dot{b}}(-g^{-1})}_{U},\]
where
\begin{align*}
&k_4=\begin{cases}
	t^du_2g^{-1}u_1^{-1}&\text{if $\dot{b}=\dot{\alpha_1}$ and $w$ is diagonal},\\
	u_2g^{-1}u_1^{-1}t^d&\text{if $\dot{b}=\dot{\alpha_1}$ and $w$ is not diagonal},\\
	u_1g^{-1}\tau^{-1}(u_2^{-1})t^d&\text{if $\dot{b}=\dot{\alpha_0}$ and $w$ is diagonal},\\
	t^d\tau(u_1g^{-1})u_2^{-1}&\text{if $\dot{b}=\dot{\alpha_0}$ and $w$ is not diagonal},
\end{cases}
\end{align*}
\end{proof}

We put $X_0=\{ (e,\tilde{w})\in E(2,\Dt)\times \tilde{N}_0\ |\ \rho_0(e)=\psi_0(\tilde{w})\}$, and define permutations $\lambda(h)$, $\mu(u)$, $\nu_{\dot{a}}$ (resp. $\lambda(h)^*$, $\mu(u)^*$, $\nu^*_{\dot{a}}$) on $X_0$ for $h\in\tilde{H}_0$, $u\in U$ and $\dot{a},\dot{b}\in\Pi_a$ as follows (see Lemma \ref{lem:212}):
\begin{align*}
\lambda(h)(e,\tilde{w})&=(\psi_0(h)e,h\tilde{w}),\\
(e,\tilde{w})\lambda(h)^*&=(e\psi_0(h),\tilde{w}h),\\
\mu(u)(e,\tilde{w})&=(ue,\tilde{w}),\\
(e,\tilde{w})\mu(u)^*&=(eu,\tilde{w}),\\
\nu_{\dot{a}}(e,\tilde{w})&=
\begin{cases}
	(w_{\dot{a}}(1)e,\tilde{w}_{\dot{a}}\tilde{w})&\text{if }\rho_0(w_{\dot{a}}(1)e)=w_{\dot{a}}(1)w,\\
	(w_{\dot{a}}(1)e,\tilde{h}_{\dot{a}}(f)^{-1}\tilde{w})&\text{if }\rho_0(w_{\dot{a}}(1)e)=h_{\dot{a}}(1)h_{\dot{a}}(f)^{-1}w,
\end{cases}\\
(e,\tilde{w})\nu^*_{\dot{b}}&=
\begin{cases}
	(ew_{\dot{b}}(-1),\tilde{w}\tilde{w}_{\dot{b}}^{-1})&\text{if }\rho_0(ew_{\dot{b}}(-1))=ww_{\dot{b}}(-1),\\
	(ew_{\dot{b}}(-1),\tilde{w}\tilde{h}_{\dot{b}}(g))&\text{if }\rho_0(ew_{\dot{b}}(-1))=wh_{\dot{b}}(g)h_{\dot{b}}(1)^{-1},
\end{cases}
\end{align*}
where
\[\tilde{w}_{\dot{a}}=
\begin{cases}
	w_1&\text{if $\dot{a}=\dot{\alpha_1}$},\\
	\tilde{h}(-t^{-1})w_1&\text{if $\dot{a}=\dot{\alpha_0}$},
\end{cases}\quad
\tilde{h}_{\dot{a}}(f)=
\begin{cases}
	\tilde{h}(f)&\text{if $\dot{a}=\dot{\alpha_1}$},\\
	\tilde{h}(-f^{-1}t^{-1})\tilde{h}(-t^{-1})^{-1}&\text{if $\dot{a}=\dot{\alpha_0}$}.
\end{cases}\]

{\rem{ In the definition of $\lambda(h)$, it is not obvious that $\rho_(\psi_0(h)e)=\psi_0(h)w$ in general. For $e=uwv\in UNU$ with some $u,v\in U$, if $\psi_0(h)e=\psi_0(h)uwv=u'\psi_0(h)wv$ by (E2), then $u'$ is not always in $U$. However, we can choose suitable $\dot{a}\in\Pi_a$, $y\in U'_{\dot{a}}$, and $x_{\dot{a}}(f)\in U_{\dot{a}}$, which satisfy $u'\in U$ in this case.
}}\vspace{3mm}

Let $G_0$ (resp. $G_0^*$) be the group of permutations on $X_0$ generated by
$\lambda(h)$, $\mu(u)$, $\nu_{\dot{a}}$ (resp. $\lambda(h)^*$, $\mu(u)^*$, $\nu^*_{\dot{a}}$) for $h\in\tilde{H}_0$, $u\in U$, and $\dot{a}\in\Pi_a$.\vspace{3mm}

{\lem\label{lem:214}{
For all $(e,\tilde{w})\in X$, $g\in G_0$, and $g^*\in G_0^*$, it holds that
\[(g(e,\tilde{w}))g^*=g((e,\tilde{w})g^*).\]
}}

\begin{proof}
It suffices to show this equality for the generators of $G_0$ and $G_0^*$. We can easily verify it, except for the case that $g=\nu_{\dot{a}}$, $g^*=\nu^*_{\dot{b}}$ for $\dot{a},\dot{b}\in\Pi_a$. It can be easily checked that the first component of $(\nu_{\dot{a}}(e,\tilde{w}))\nu^*_{\dot{b}}$ is equal the one of $\nu_{\dot{a}}((e,\tilde{w})\nu^*_{\dot{b}})$. Let us show that the second components of them are equal. We write $e=yx_{\dot{a}}(-f)wx_{\dot{b}}(g)z$ with some $f,g\in D$, $w\in N$, $y\in U'_{\dot{a}}$, $z\in U'_{\dot{b}}$, and $\dot{a},\dot{b}\in\Pi_a$, as in Lemma \ref{lem:212}. Here, we recall from ($\spadesuit$) that $w$ is one of the following two cases; for $u_i\in D_{\tau}^{\times}$,
\begin{equation*}
\underbrace{w=\left(\begin{array}{cc}
u_1&0\\
0&u_2
\end{array}\right)}_{\text{diagonal}}\quad\text{or}\quad
\underbrace{w=\left(\begin{array}{cc}
0&u_2\\
u_1&0
\end{array}\right)}_{\text{anti-diagonal}},
\end{equation*}
and that $d=d_w=\mathrm{deg}(u_1^{-1}u_2)$. We denote the second component of $(\nu_{\dot{a}}(e,\tilde{w}))\nu^*_{\dot{b}}$ (resp. $\nu_{\dot{a}}((e,\tilde{w})\nu^*_{\dot{b}})$) by $(\nu_{\dot{a}}(\tilde{w}))\nu^*_{\dot{b}}$ (resp. $\nu_{\dot{a}}((e,\tilde{w})\nu^*_{\dot{b}})$) for convenience. We note that the following relation holds for all $u\in\Dt^{\times}$.
\[(\star)\hspace{3mm}w_1\tilde{h}(u)w_1^{-1}=\tilde{h}(u^{-1})=\tilde{h}(-u)^{-1}\tilde{h}(-1).\]

\fbox{Case 1} Assume that $w(\dot{b})\neq \pm \dot{a}$.\vspace{3mm}

If we put $\Gamma(\sigma_{\dot{c}})=\{\dot{d}\in\Delta_a^+\ |\ \sigma_{\dot{c}}(\dot{d})\in\Delta_a^-\}$ for $\dot{c}\in\Pi_a$, then $\Gamma(\sigma_{\dot{c}})=\{\dot{c}\}$. Therefore, we have the following trivial cases:
\begin{align*}
&(\nu_{\dot{a}}(\tilde{w}))\nu^*_{\dot{b}}=\nu_{\dot{a}}((\tilde{w})\nu^*_{\dot{b}})\\
&=\begin{cases}
\tilde{w}_{\dot{a}}\tilde{w}\tilde{w}_{\dot{b}}^{-1}&\text{if ``$w^{-1}(\dot{a})\in\Delta_a^+$ or $f=0$'' and ``$w(\dot{b})\in\Delta_a^+$ or $g=0$''},\\
\tilde{w}_{\dot{a}}\tilde{w}\tilde{h}_{\dot{b}}(g)&\text{if ``$w^{-1}(\dot{a})\in\Delta_a^+$ or $f=0$'' and ``$w(\dot{b})\not\in\Delta_a^+$''},\\
\tilde{h}_{\dot{a}}(f)^{-1}\tilde{w}\tilde{w}_{\dot{b}}^{-1}&\text{if ``$w^{-1}(\dot{a})\not\in\Delta_a^+$'' and ``$w(\dot{b})\in\Delta_a^+$ or $g=0$''},\\
\tilde{h}_{\dot{a}}(f)^{-1}\tilde{w}\tilde{h}_{\dot{b}}(g)&\text{if ``$w^{-1}(\dot{a})\not\in\Delta_a^+$'' and ``$w(\dot{b})\not\in\Delta_a^+$''}.
\end{cases}
\end{align*}

\fbox{Case 2} Assume that $w(\dot{b})=\dot{a}$.\vspace{3mm}

We see that $\dot{a}=\dot{b}$, $d=0$, and $w$ is diagonal. It suffices to show the claim when $\dot{a}=\dot{\alpha_0}$. If $f=0$ or $g=0$, then we see that the claim follows from Lemma \ref{lem:26} and Lemma \ref{lem:210}. We assume that $f\neq 0$ and $g\neq 0$.
\begin{align*}
w^{-1}x_{\dot{\alpha_0}}(-f)w&=x_{w^{-1}(\dot{\alpha_0})}(-\tau^{-1}(u_2^{-1})fu_1),\\
wx_{\dot{\alpha_0}}(g)w^{-1}&=x_{w(\dot{\alpha_0})}(\tau^{-1}(u_2)gu_1^{-1}).
\end{align*}
We put $s=-f+\tau^{-1}(u_2)gu_1^{-1}$ and $s'=-\tau^{-1}(u_2^{-1}fu_1)+g$. If $s=0$, then $s'=0$, and the claim is obvious. We assume $s\neq 0$. We have to show that
\[\tilde{h}_{\dot{\alpha_0}}(-s)^{-1}\tilde{w}\tilde{w}_{\dot{\alpha_0}}^{-1}=\tilde{w}_{\dot{\alpha_0}}\tilde{w}\tilde{h}_{\dot{\alpha_0}}(s'),\]
which is equivalent to
\[\tilde{h}(s^{-1}t^{-1})\tilde{w}w_1=w_1\tilde{w}\tilde{h}(-u_1^{-1}s^{-1}t^{-1}u_2).\]
We compute
\begin{align*}
w_1\tilde{w}\tilde{h}&(-u_1^{-1}s^{-1}t^{-1}u_2)\\
&=w_1\tilde{h}(-s^{-1}t^{-1})\tilde{h}(u_1u_2^{-1})^{-1}\tilde{w}&\text{(by Lemma \ref{lem:26})}\\
&=\tilde{h}(s^{-1}t^{-1})^{-1}\tilde{h}(-u_1u_2^{-1})w_1\tilde{w}&\text{(by ($\star$))}\\
&=\tilde{h}(s^{-1}t^{-1})^{-1}\tilde{h}(-u_1u_2^{-1})\tilde{h}(u_1u_2^{-1})^{-1}\tilde{w}w_1&\text{(by Lemma \ref{lem:210})}\\
&=\tilde{h}(-s^{-1}t^{-1})^{-1}\tilde{h}(-u_1u_2^{-1})\tilde{h}(-u_1u_2^{-1})^{-1}\tilde{w}w_1^{-1}&\text{(by ($\star$))}\\
&=\tilde{h}(-s^{-1}t^{-1})^{-1}\tilde{w}w_1^{-1}\\
&=\tilde{h}(s^{-1}t^{-1})\tilde{w}w_1&\text{(by ($\star$))}.
\end{align*}

\fbox{Case 3} Assume that $w(\dot{b})=-\dot{a}$.\vspace{3mm}

We see that $\dot{a}=\dot{b}$, $d=-2$, and $w$ is not diagonal. It suffices to show the claim when $\dot{a}=\dot{\alpha_0}$. If $f=0$ or $g=0$, then we see that the claim follows from Lemma \ref{lem:26} and Lemma \ref{lem:210}. We assume $f\neq 0$ and $g\neq 0$. 
\begin{align*}
w^{-1}x_{\dot{\alpha_0}}(-f)w&=x_{w^{-1}(\dot{\alpha_0})}(-\tau^{-1}(u_2^{-1})fu_1),\\
wx_{\dot{\alpha_0}}(g)w^{-1}&=x_{w(\dot{\alpha_0})}(\tau^{-1}(u_2)gu_1^{-1}).
\end{align*}
Put $s=-t^{-1}u_2^{-1}f^{-1}t^{-1}u_1+g$ and $s'=f-t^{-1}u_1g^{-1}t^{-1}u_2^{-1}$. If $s=0$, then $s'=0$, and the claim follows from Case 2. We assume that $s\neq 0$. We have to show that
\[\tilde{h}_{\dot{\alpha_0}}(f)^{-1}\tilde{w}\tilde{h}_{\dot{\alpha_0}}(s)
=\tilde{h}_{\dot{\alpha_0}}(s')^{-1}\tilde{w}\tilde{h}_{\dot{\alpha_0}}(g),\]
but this can be written as
\[\tilde{h}(-f^{-1}t^{-1})^{-1}\tilde{w}\tilde{h}(-s^{-1}t^{-1})=\tilde{h}(-s'^{-1}t^{-1})^{-1}\tilde{w}\tilde{h}(-g^{-1}t^{-1}).\]
We see that
\begin{align*}
&\tilde{h}(-f^{-1}t^{-1})^{-1}\tilde{w}\tilde{h}(-s^{-1}t^{-1})=\tilde{h}(-s'^{-1}t^{-1})^{-1}\tilde{w}\tilde{h}(-g^{-1}t^{-1})\\
&\Longleftrightarrow\tilde{w}\tilde{h}(-s^{-1})\tilde{h}(-g^{-1}t^{-1})^{-1}\tilde{w}^{-1}=\tilde{h}(-f'^{-1}t^{-1})\tilde{h}(-s'^{-1}t^{-1})^{-1}\\
&\Longleftrightarrow (w_1\tilde{w})\tilde{h}(-s^{-1})\tilde{h}(-g^{-1}t^{-1})^{-1}(w_1\tilde{w})^{-1}=w_1\tilde{h}(-f'^{-1}t^{-1})\tilde{h}(-s'^{-1}t^{-1})^{-1}w_1^{-1}\\
&\Longleftrightarrow\tilde{h}(u_1s^{-1}t^{-1}u_2^{-1})\tilde{h}(u_1g^{-1}t^{-1}u_2^{-1})^{-1}=\tilde{h}(f^{-1}t^{-1})^{-1}\tilde{h}(s'^{-1}t^{-1})\\
&\hspace{7cm}\text{(by Lemma \ref{lem:26} and ($\star$))}\\
&\Longleftrightarrow\tilde{h}(f^{-1}t^{-1})\tilde{h}(u_1s^{-1}t^{-1}u_2^{-1})=\tilde{h}(s'^{-1}t^{-1})\tilde{h}(u_1g^{-1}t^{-1}u_2^{-1})\\
&\Longleftrightarrow\tilde{h}(tf)\tilde{h}(u_2tsu_1^{-1})=\tilde{h}(ts')\tilde{h}(u_2tgu_1^{-1})\hspace{3cm}\text{(by ($\star$))}
\end{align*}
If we put $x=u_2tgu_1^{-1}$ and $y=tf$, then the last identity is equivalent to
\[c(y,x-y^{-1})=c(y-x^{-1},x),\]
and we obtain by $\mathrm{deg}(xy)=0$
\begin{align*}
c(y,x-y^{-1})&=c(y,(xy-1)y^{-1})\\
	&=c(y(xy-1),y^{-1})&(\text{by (P2)})\\
	&=c((y-x^{-1})xy,y^{-1})\\
	&=c(-(y-x^{-1})x,y^{-1})&(\text{by (P5)$'$})\\
	&=c(1-yx,y^{-1})\\
	&=c(1-yx,x)&(\text{by (P4)})\\
	&=c(-(y-x^{-1})x,x)\\
	&=c(y-x^{-1},x).&(\text{by (P5)$'$})
\end{align*}
Thus we have proved Lemma \ref{lem:214}.
\end{proof}\vspace{3mm}

{\lem\label{lem:215}{ The actions of $G_0$ and $G_0^*$ are simply transitive on $X_0$.}}\vspace{3mm}

\begin{proof}
We give the proof only for $G_0$ because the one for $G_0^*$ is similar. We first show the transitivity. Let $(e,\tilde{w}),(e',\tilde{w}')\in X$. Since $E(2,\Dt)$ is generated by $U$ and $w_{12}(1)$, we deduce that there exists $g'\in\langle \mu(u),\nu_{\dot{a}}\ |\ u\in U,\ \dot{a}\in\Pi_a\rangle\subset G_0$ such that $g'(e,\tilde{w})=(e',\tilde{w}^*)$ for some $\tilde{w}^*$. $(e',\tilde{w}')$ and $(e',\tilde{w}^*)$ lie in $X_0$, so we see that $\tilde{w}^*= l\tilde{w}$ for some $l\in L_0$, that is, $(e',\tilde{w}^*)=\lambda(l)(e',\tilde{w}')$. Thus, there exists $g\in G_0$ such that $g(e,\tilde{w})=(e',\tilde{w}')$.

If there are $g_1,g_2\in G_0$ with $g_1x=g_2x$ for some $x\in X_0$, then we have $g_1(xg^*)=g_2(xg^*)$ for any $g^*\in G_0^*$ by Lemma 2.14. This implies that $g_1x'=g_2x'$ for every $x'\in X_0$ by the transitivity of $G_0^*$, which yields $g_1=g_2$.
\end{proof}

{\thm\label{thm:216}{ There exists a surjective group homomorphism $\Psi_0:G_0\to E(2,\Dt)$ sending $\lambda(h)$ to $\psi_0(h)$ for $h\in\tilde{H}_0$, $\mu(u)$ to $u$ for $U$ and $\nu_{\dot{\alpha}}$ to $w_{\dot{\alpha}}(1)$ for $\dot{\alpha}\in\Pi_a$, and which satisfies that the following exact sequence is a central extension of $E(2,\Dt)$.
\[1\longrightarrow L_0\longrightarrow G_0\stackrel{\Psi_0}{\longrightarrow} E(2,\Dt)\longrightarrow 1\]}}
\begin{proof}
We only prove that $\Psi_0$ is a central extension since the first assertion is obvious from Lemma \ref{lem:215}. If $g\in\mathrm{Ker}\ \Psi_0$, then $g(e,\tilde{w})=(e,\tilde{w}')$. $\rho_0(e)=\psi_0(\tilde{w})=\psi_0(\tilde{w}')$ implies $\tilde{w}'=l\tilde{w}$ for some $l\in L_0$, that is, $g(e,\tilde{w})=\lambda(l)(e,\tilde{w})$. We obtain $g=\lambda(l)$ by Lemma \ref{lem:215}. Therefore $\mathrm{Ker}\ \Psi_0\cong L_0$. Let $(e,\tilde{w})\in X_0$. For $(e,\tilde{w}), (e',\tilde{w}')\in X_0$ and $g\in G_0$, if $g(e,\tilde{w})=(e',\tilde{w}')$, then we see that $\lambda(l)g(e,\tilde{w})=(e',l\tilde{w}')=g\lambda(l)(e,\tilde{w})$.
\end{proof}\vspace{3mm}

\begin{proof}[PROOF OF THEOREM \ref{thm:23}]
By Theorem \ref{thm:216} and \cite[Theorem 3]{rs}, we get the following diagram.
\[\begin{split}
\xymatrix{
&1 \ar[r] &L_0 \ar[r]\ar[d]^-{\zeta_0} &G_0 \ar[r]^-{\Psi_0}  &E(2,D_{\tau})\ar@{=}[d] \ar[r] &1\\
&1 \ar[r]&K_2(2,\Dt) \ar[r] &St(2,D_{\tau}) \ar[r]\ar[u]^-{p_0}  &E(2,D_{\tau}) \ar[r] &1
}
\end{split}\]
Then, we see in a way similar way to \cite[Proposition 4.9]{r2} that there exists a unique surjective group homomorphism $p_0: St(2,\Dt)\to G_0$, so that we find that $p_0(K_2(2,\Dt))\subset L_0$ and $p_0(\hat{c}(u,v))=c(u,v)$ for all $u,v\in\Dt^{\times}$. On the other hand, we know $\zeta_0$ sends $c(u,v)$ to $\hat{c}(u,v)$ for $u,v\in D_{\tau}$. If we consider the generators of $L_0$ and $K_2(2,\Dt)$, then we obtain $L_0\cong K_2(2,\Dt)$.
\end{proof}

{\cor{ We have the following exact sequence:
\[1\longrightarrow K_2(2,\Dt)\longrightarrow P\stackrel{\varphi_0}{\longrightarrow} \Dt^{\times}\longrightarrow K_1(2,\Dt)\longrightarrow 1.\]}}


\section{Non-symplectic $K_2$-groups}

\hspace{5mm}Let $D_{\tau}$ be the ring of non-commutative Laurent polynomials defined in Section 1. Let $Q$ be the group generated by $c(u,v)$, $u,v\in D_{\tau}^{\times}$ with the defining relations that for $u,v,w\in\Dtt$ and $s\in D^{\times}$ with $s\neq 1$,
\begin{align*}
&\text{(Q1)}\quad c(uv,w)=c(^uv,^uw)c(u,w),\\
&\text{(Q2)}\quad c(u,vw)=c(u,v)c(^vu,^vw),\\
&\text{(Q3)}\quad c(s,1-s)=1,
\end{align*}
whiere we write $^uv=uvu^{-1}$. We also denote by $^xc(u,v)=c(^xu,^xv)$ for $x,u,v\in\Dtt$. Since $[u,v]\in [D_{\tau}^{\times},D_{\tau}^{\times}]$, $u,v\in D_{\tau}^{\times}$, satisfy the same relations (Q1)--(Q3) (with $c(u,v)$ replaced by $[u,v]$), there exists a (unique) surjective group homomorphism $\varphi: Q\twoheadrightarrow [D_{\tau}^{\times},D_{\tau}^{\times}]$ which sends $c(u,v)$ to $[u,v]$ for $u,v\in D_{\tau}^{\times}$. Set $L=\mathrm{Ker}\ \varphi$; note that
\begin{align*}
L=\{ c(u_1,v_1)^{p_1}c(u_2,v_2)^{p_2}\cdots &c(u_r,v_r)^{p_r}\ |\ p_i=\pm 1,\ u_i,v_i\in\Dt^{\times},\\
&[u_1,v_1]^{p_1}[u_2,v_2]^{p_2}\cdots [u_r,v_r]^{p_r}=1\}.&\text{(\#\#)}
\end{align*}

{\prop\label{prop:31}{ The following exact sequence is a central extension of $Q$ by $L$.
\[1\longrightarrow L\longrightarrow Q\stackrel{\varphi}{\longrightarrow} [D_{\tau}^{\times},D_{\tau}^{\times}]\longrightarrow 1.\]}}

\begin{proof}
First, we prove the following relation:
\[\text{(Q4)}\quad c(u,v)c(x,y)=^{[u,v]}c(x,y)c(u,v)\]
Let $a,b,u,v\in\Dt^{\times}$. Using (Q1),
\begin{align*}
c(ua,vb)&=c(^ua,^uvb)c(u,vb)\\
		&=c(^ua,^ub)c(^{uv}a,^{uv}b)c(u,v)c(^vu,^vb).
\end{align*}
On the other hand, by (Q2),
\begin{align*}
c(ua,vb)&=c(ua,v)c(^vua,^vb)\\
		&=c(^ua,^uv)c(u,v)c(^{vu}a,^{vu}b)c(^vu,^vb).
\end{align*}
Therefore, we get $c(u,v)c(^{vu}a,^{vu}b)=c(^{uv}a,^{uv}b)c(u,v)$. Putting $x= ^{vu}\! a$ and $y= ^{vu}\! b$ respectively, we see that $^{uv}a=^{[u,v]}x$, $^{uv}b=^{[u,v]}y$. Hence (Q4) holds.\vspace{3mm}

Applying (Q4) twice, we obtain
\[c(u_2,v_2)c(u_1,v_1)c(x,y)c(u_1,v_1)^{-1}c(u_2,v_2)^{-1}=^{[u_2,v_2][u_1,v_1]}c(x,y).\]
If we take $\xi=c(u_1,v_1)^{p_1}c(u_2,v_2)^{p_2}\cdots c(u_r,v_r)^{p_r}\in Q$ for $u_i,v_i\in\Dtt$ and $p_i\in\{\pm 1\}$, using an induction on $r$, we see that
\[\xi c(u,v)\xi^{-1}=\ ^{\varphi(\xi)}c(u,v).\]
Thus $\xi$ is central if $\varphi(\xi)=1$.
\end{proof}\vspace{3mm}

By comparing (Q1), (Q2), (Q3) with (TT6), (TT7), (TT3) recpectively, we see that there exists a (unique) surjective group homomorphism $\zeta: Q\to K_2(n,D_{\tau})$ with maps $c(u,v)$ to $\hat{c}(u,v)$ for $u,v\in D_{\tau}^{\times}$. By Proposition \ref{prop:31} and ($\#\#$), the restriction of this $\zeta$ to $L\subset Q$ is a surjective group homomorphism from $L$ onto $K_2(n,D_{\tau})$. \vspace{3mm}

{\thm\label{thm:32}{
The group homomorphism $\zeta: L \rightarrow K_{2}(n,D_{\tau})$ is
an isomorphism of groups.
}}\vspace{3mm}

Let $\tilde{H}$ be the group presented by $\tilde{h}_{ij}(u)$ for $u\in\Dt^{\times}$ and $1\leq i\neq j\leq n$ and $z(q)$, $q\in Q$ with the following relations that
\begin{align*}
\text{(H1)}\quad &\tilde{h}_{ij}(u)\tilde{h}_{ji}(u)=1,\\
\text{(H2)}\quad &\tilde{h}_{ij}(u)\tilde{h}_{ki}(u)\tilde{h}_{jk}(u)=1,\\
\text{(H3)}\quad &\tilde{h}_{ij}(u)\tilde{h}_{ik}(v)\tilde{h}_{ij}(u)^{-1}=\tilde{h}_{ik}(uv)\tilde{h}_{ik}(u)^{-1}&(k\neq j),\\
\text{(H4)}\quad &\tilde{h}_{ij}(u)\tilde{h}_{kj}(v)\tilde{h}_{ij}(u)^{-1}=\tilde{h}_{kj}(vu)\tilde{h}_{kj}(u)^{-1}&(k\neq i),\\
\text{(H5)}\quad &[\tilde{h}_{ij}(u),\tilde{h}_{kl}(v)]=1&(k\neq i,j\text{ and } l\neq i,j),\\
\text{(H6)}\quad &z(q)z(q')=z(qq'),\\
\text{(H7)}\quad &\tilde{h}_{ij}(u)\tilde{h}_{ij}(v)\tilde{h}_{ij}(vu)^{-1}=z(c(u,v)),\\
\text{(H8)}\quad &\tilde{h}_{ij}(u)\tilde{h}_{ik}(v)\tilde{h}_{ij}(u)^{-1}\tilde{h}_{ik}(v)^{-1}=z(c(u,v)),\\
\text{(H9)}\quad &\tilde{h}_{ij}(u)z(c(u,v))=z(^xc(u,v))\tilde{h}_{ij}(u)
\end{align*}
for $1\leq i\neq j\leq n$, $1\leq k\neq l\leq n$, and $u,v\in\Dtt$. We deduce by (H6) that $\{z(q)\ |\ q\in Q\}$ is a subgroup of $\tilde{H}$. Moreover, we see by a similar to \cite[CHAPITRE II]{hm} and \cite[Proposition 2]{jmur} that $L$ is isomophic to $\{z(l)\ |\ l\in L\}$. We identify $\{z(l)\ |\ l\in L\}$ with $L$, and write $z(l)$ simply by $l$ for  $l\in L$.\vspace{3mm}

For $u,v\in\Dtt$ and $1\leq i\neq j\leq n$, we find that
\begin{align*}
\tilde{h}_{ij}(u)&\tilde{h}_{ij}(v)\tilde{h}_{ij}(u)^{-1}\\
&=\tilde{h}_{kj}(u)\tilde{h}_{ik}(u)\tilde{h}_{ij}(v)\tilde{h}_{ik}(u)^{-1}\tilde{h}_{kj}(u)^{-1}\quad&\text{(by (H2))}\\
&=\tilde{h}_{ij}(uvu)\tilde{h}_{ij}(u^2)^{-1}\quad&\text{(by (H3) and (H4))}
\end{align*}
for $1\leq k\leq n$ with $k\neq i$ and $k\neq j$. If we put $v=u^{-1}$, then
\[\tilde{h}_{ij}(u^{-1})\tilde{h}_{ij}(u)^{-1}=\tilde{h}_{ij}(u^2)^{-1}.\]
Therefore, we get $\tilde{h}_{ij}(u)\tilde{h}_{ij}(v)=\tilde{h}_{ij}(uvu)\tilde{h}_{ij}(u^{-1})$ for $1\leq i\neq j\leq n$ and $u,v\in\Dtt$. Using this relation, we can easily check the following proposition.\vspace{3mm}

{\prop\label{prop:33}{
All relations in Lemma \ref{lem:14} and Lemma \ref{lem:15} with $\hat{h}_{ij}(u)$ replaced by $\tilde{h}_{ij}(u)$ and $\hat{c}_{ij}(u,v)$ replaced by $z(c(u,v))$ for $u,v\in\Dtt$ hold in $\tilde{H}$.
}}\vspace{3mm}

{\prop\label{prop:34}{There exists a (unique) surjective group homomorphism $\pi: \tilde{H}\to T$ which sends $\tilde{h}_{ij}(u)$ to $h_{ij}(u)$ for $u\in D_{\tau}^{\times}$. The kernel $\mathrm{Ker}\ \pi$ is identical to $L$. Moreover, 
\[1\longrightarrow L\longrightarrow \tilde{H}\stackrel{\pi}{\longrightarrow} T\longrightarrow 1\]
 is a central extension of $T$ by $L$.}}\vspace{3mm}

\begin{proof}
We first note that the following relations hold: for $1\leq i\neq j\leq n$ and $1\leq k\leq n$ with $k\neq i$, $k\neq j$,
\begin{align*}
\tilde{h}_{ij}(u)\tilde{h}_{ij}(v)&\equiv\tilde{h}_{ij}(vu)&\text{mod $z(Q)$},\\
\tilde{h}_{ij}(u)\tilde{h}_{ik}(v)&\equiv\tilde{h}_{ik}(v)\tilde{h}_{ij}(u)&\text{mod $z(Q)$},
\end{align*}
where $z(Q)=\{z(q)\ |\ q\in Q\}$. Therefore, if we take $\tilde{h}\in\mathrm{Ker}\ \pi$, then $\tilde{h}$ is expressed as
\[\tilde{h}\equiv\tilde{h}_{12}(u_1)\tilde{h}_{13}(u_2)\cdots\tilde{h}_{1n}(u_{n-1})\hspace{5mm}\text{mod $z(Q)$}\]
for $u_i\in\Dtt$. Applying $\pi$ to $\tilde{h}$, we have $\pi(\tilde{h})=\mathrm{diag}(u_1\cdots u_{n-1},u_1^{-1},\dots,u_{n-1}^{-1})$, but this implies $u_i=1$ for all $1\leq i\leq n-1$. Hence $\tilde{h}\equiv 1$ (mod $z(Q)$). Thus $\tilde{h}\in\mathrm{Ker}\ \pi$ is of the form
\[\tilde{h}=z(c(u_1,v_1)c(u_2,v_2)\cdots c(u_r,v_r))\in z(Q)\]
for $u_i,v_i\in\Dtt$, and we see that $\tilde{h}\in\mathrm{Ker}\ \pi\subset L$. It is obvious that $L\subset \mathrm{Ker}\ \pi$.
\end{proof}\vspace{3mm}

\noindent By Propositions \ref{prop:31} and \ref{prop:34}, we get the commutative diagram;
\begin{equation}
\begin{split}
\xymatrix{
&L \ar[r]\ar@{=}[d] &Q \ar[r]^-{\varphi}\ar[d]^-{z}  &[D_{\tau}^{\times},D_{\tau}^{\times}]\ar[d]^-{d}\\
&L \ar[r] &\tilde{H} \ar[r]^-{\pi}  &T
}
\end{split}
\tag{CD2}
\end{equation}
where $d$ is an embedding defined by $d([u,v])=\mathrm{diag}([u,v],1,\dots,1)$. This implies that $Q$ is isomorphic to $\{z(q)\ |\ q\in Q\}$. We identify $\{z(q)\ |\ q\in Q\}$ with $Q$, and write $z(q)$ simply by $q$ for $q\in Q$.\vspace{3mm}

Next, we construct some extension of the monomial subgroup $N$, which is ``compatible'' with the extension $(\tilde{H},\pi)$ of $T$ in Proposition \ref{prop:34} (see Proposition \ref{prop:37} below). For this, we give the presentation of $N$, and then construct an action of $N$ on $\tilde{H}$. The next lemma and proposition follow from \cite[Proposition 3]{jmur} and \cite[Proposition 5.8]{hs}.\vspace{3mm}

{\lem\label{lem:35}{
The subgroup $N$ of $E(n,\Dt)$ is the group generated by $w_{\alpha}(u)$ for $u\in \Dtt$ and $\alpha\in\Pi$, with the following defining relations that for $\alpha=\epsilon_i-\epsilon_j,\beta=\epsilon_k-\epsilon_l\in\Pi$ and $u,v\in\Dtt$,}}
\begin{align*}
&\text{(N1)}\hspace{5mm} w_{\alpha}(u)^{-1}=w_{\alpha}(-u),\\
&\text{(N2)}\hspace{5mm} w_{\alpha}(1)h_{\beta}(u)w_{\alpha}(1)^{-1}=h_{\beta}(u)h_{\alpha}(u^{-\langle\alpha,\beta\rangle}),\\
&\text{(N3)}\hspace{5mm} h_{\alpha}(u)h_{\alpha}(v)=h_{\alpha}(uvu)h_{\alpha}(u^{-1}),\\
&\text{(N4)}\hspace{5mm} h_{\alpha}(u)h_{\beta}(v)h_{\alpha}(u)^{-1}=
	\begin{cases}
	h_{\beta}(u^{-1}v)h_{\beta}(u)&\text{if $j=k$},\\
	h_{\beta}(vu^{-1})h_{\beta}(u)&\text{if $i=l$},\\
	h_{\beta}(v)&\text{if $\alpha+\beta\not\in\Delta$},
	\end{cases}\\
&\text{(N5)}\hspace{5mm} w_{\alpha}(1)w_{\beta}(1)w_{\alpha}(1)=w_{\beta}(1)w_{\alpha}(1)w_{\beta}(1),\hspace{20mm}(\langle\alpha,\beta\rangle=-1)\\
&\text{(N6)}\hspace{5mm} w_{\alpha}(1)w_{\beta}(1)=w_{\beta}(1)w_{\alpha}(1).\hspace{38mm}(\langle\alpha,\beta\rangle=0)
\end{align*}

{\prop\label{prop:36}{
There exists an action of $N$ on $\tilde{H}$ defined by
\begin{align*}
w_{ij}(u)&\cdot \tih_{kl}(v)\\
&=\begin{cases}
\tih_{kl}(v)&\text{if $i\neq k,l$ and $j\neq k,l$}\\
\tih_{ji}(-u^{-1}vu^{-1})\tih_{ji}(-u^{-2})^{-1}&\text{if $i=k$ and $j=l$}\\
\tih_{ij}(-uvu)\tih_{ji}(-u^2)^{-1}&\text{if $i=l$ and $j=k$}\\
\tih_{jl}(-u^{-1}v)\tih_{jl}(-u^{-1})^{-1}&\text{if $i=k$ and $j\neq l$}\\
\tih_{kj}(-vu)\tih_{kj}(-u)^{-1}&\text{if $i=l$ and $j\neq k$}\\
\tih_{il}(uv)\tih_{il}(u)^{-1}&\text{if $i\neq l$ and $j=k$}\\
\tih_{ki}(vu^{-1})\tih_{ki}(u^{-1})^{-1}&\text{if $i\neq k$ and $j=l$}
\end{cases}
\end{align*}
for $1\leq i\neq j\leq n$, $1\leq k\neq l\leq n$, and $u,v\in\Dt^{\times}$. We also denote $w_{ij}(u)\cdot h_{kl}(v)$ by $w_{ij}(u)\tih_{kl}(v)w_{ij}(u)^{-1}$ for convenience.
}}\vspace{3mm}

Let $\tilde{W}$ be the the group generated by $\tiw_{\alpha}$ for all $\alpha\in \Pi$ with the following defining relations that
\begin{align*}
&\text{(W1)}\quad \tih_{\alpha}\tiw_{\beta}\tih_{\alpha}^{-1}=\tiw_{\beta}^{\delta},\\
&\text{(W2)}\quad \tiw_{\alpha}\tiw_{\gamma}\tiw_{\alpha}=\tiw_{\gamma}\tiw_{\alpha}\tiw_{\gamma}&\text{if $\langle\alpha,\gamma\rangle=-1$},\\
&\text{(W3)}\quad \tiw_{\alpha}\tiw_{\gamma}=\tiw_{\gamma}\tiw_{\alpha}&\text{if $\langle\alpha,\gamma\rangle=0$},
\end{align*}
where $\alpha\neq\pm\gamma$, $\delta=(-1)^{\langle\alpha,\beta\rangle}$, and $\tih_{\alpha}=\tiw_{\alpha}^2$. Put $\tilde{T}=\langle \tih_{\alpha}\ |\ \alpha\in\Pi\rangle$ and $\tilde{N}^*=\tilde{H}\rtimes\tilde{W}$, where $\tilde{W}$ acts on $\tilde{H}$ by $\tiw_{\alpha}\cdot\tih=w_{\alpha}(-1)\cdot\tih$ for $\alpha\in\Pi$ and $\tih\in\tilde{H}$. Then we see by Proposition \ref{prop:36} that $\tilde{T}$ is the group generated by $\tih_{\alpha}$ for all $\alpha\in\Pi$ with the defining relation: (T)\quad $\tih_{\alpha}\tih_{\beta}\tih_{\alpha}^{-1}=\tih_{\beta}^{\delta}$. Then there is a canonical homomorphism $\eta:\tilde{T}\to\tilde{H}$ defined by $\eta(\tih_{\alpha})=\tih_{\alpha}(-1)$.  Let $J$ be the subgroup of $\tilde{N}^*$ generated by $(\tilde{t},\eta(\tilde{t})^{-1})$ for all $\tilde{t}\in\tilde{T}$, and let $\tilde{N}=\tilde{N}^*/J$ be the quotient group. Let $\tiw_{\alpha}J$ be the canonical image of $\tiw_{\alpha}$ in $\tilde{N}$. Let $\tilde{\psi}:N^*\to \tilde{N}$ be the canonical homomorphism. We deduce that the restriction of $\tilde{\psi}$ to $\tilde{H}$ is injective; we regard $\tilde{H}$ as a subgroup of $\tilde{N}$. If we set $\tiw_{\alpha}(u)=\tih_{\alpha}(u)\tiw_{\alpha}^{-1}J\in\tilde{N}$ for $u\in\Dtt$, then we find that $\tiw_{\alpha}(-1)=\tiw_{\alpha}$, $\tiw_{\alpha}(u)^{-1}=\tiw_{\alpha}(-u)$ and $\tih_{\alpha}(u)=\tiw_{\alpha}(u)\tiw_{\alpha}(-1)$ hold in $\tilde{N}$ (see Section 2). Notice that there exists a group homomorphism $\psi^*$ from $\tilde{N}^*$ to $N$ such that $\psi^*(\tilde{w}_{\alpha})=w_{\alpha}(-1)$ for $\alpha\in\Pi$ and $\psi^*(\tilde{h})=\pi(\tilde{h})$ for $\tilde{h}\in\tilde{H}$ since it can easily checked that relations in $\tilde{N}^*$ hold in $N$. We see that $J\subset\mathrm{Ker}\ \psi^*$, and hence, $\psi^*$ induces a homomorphism $\psi:\tilde{N}\to N$ which sends $\tiw_{\alpha}(u)$ to $w_{\alpha}(u)$ for $u\in\Dtt$. We see that $\psi$ is surjective by (R6), and that $\psi$ is a central extension of $N$ by $L$. Using this and Proposition \ref{prop:34}, we deduce that the restriction of $\psi$ to $\tilde{H}$ is injective (see (CD2)). Thus, we obtain the following proposition..\vspace{3mm}

{\prop\label{prop:37}{ The kernel $\mathrm{Ker}\ \psi$ of the group homomorphism $\psi: \tilde{N}\to N$ is contained in the center of $\tilde{N}$, and is isomorphic to $L$. Namely,
\[1\longrightarrow L\longrightarrow \tilde{N}\stackrel{\psi}{\longrightarrow} N\longrightarrow 1\]
is a central extension of $N$ by $L$. Moreover, the restriction of $\psi$ to $\tilde{H}$ coincides with the group homomorphism $\pi:\tilde{H}\to T\subset N$ defined in Proposition \ref{prop:34}.}}\vspace{3mm}

{\lem[cf. Lemma \ref{lem:211}]\label{lem:38}{ Every matrix $e\in E(n,\Dt)$ can be written as $e=uwv$ with some $u,v\in U$ and $w\in N$. Moreover, the monomial matrix part $w$ is uniquely determined by $e$; we define $\rho: E(n,\Dt)\to N$ by $\rho(e)=\rho(uwv)=w$.}\vspace{3mm}

Next, we determine the value of $\rho$ after multiplying a double coset by $w_{\dot{\alpha}}(\pm 1)$ for all $\dot{\alpha}\in\Pi_a$. For all $w\in N$, we can express $w=P_{\sigma}\mathrm{diag}(u_1,\dots,u_n)$ with suitable $u_i\in\Dtt$, where $P_{\sigma}$ is the permutation matrix corresponding to some permutation $\sigma$ of $\{1,\dots,n\}$. In what follows, we write $\dot{\beta}=(ij,m)$ for $\dot{\beta}=(\beta,m)\in\Delta_a$ with $\beta=\epsilon_i-\epsilon_j\in\Delta$. Then, we see by (R3) and (R4) that for $\dot{\beta}=(\beta,m)$, $\beta=\epsilon_i-\epsilon_j$,
\begin{align*}
w^{\pm 1}&x_{\dot{\beta}}(f)w^{\mp 1}\\
	&=w^{\pm 1}x_{(ij,m)}(f)w^{\mp 1}\\
	&=\begin{cases}
		x_{(\sigma(ij),m\mp d)}(\tau^{m\mp d}(t^{\pm d}\tau^{-m}(u_i^{\pm 1}f)u_j^{\mp 1}))&\text{if }\beta\in\Delta^+,\\
		x_{(\sigma^{-1}(ij),m\mp d')}(\tau^{-m\pm d'}(u_{\sigma^{-1}(i)}^{\pm 1}\tau^m(fu_{\sigma^{-1}(j)}^{\mp 1})t^{\pm d'}))&\text{if }\beta\in\Delta^-,
	\end{cases}
\end{align*}
where $d=\text{deg}(u_iu_j^{-1})$, $d'=\text{deg}(u_{\sigma^{-1}(i)}u_{\sigma^{-1}(j)}^{-1})$, and $\sigma(ij)=\sigma(i)\sigma(j)$. For $f\in D$ and $\dot{\beta}\in\Delta_a$, if $wx_{\dot{\beta}}(f)w^{-1}=x_{\dot{\gamma}}(g)$ for suitable $g\in D$ and $\dot{\gamma}\in\Delta_a$, then we denote $\dot{\gamma}$ by $w(\dot{\beta})$. Here, we know from Proposition \ref{prop:11} and \cite[Proposition 2]{rs} that every element $e\in E(n,D_{\tau})$ can be written in the form $e=yx_{{\dot{a}}}(-f)wx_{{\dot{b}}}(g)z$ with some $f,g\in D$, ${\dot{a}},{\dot{b}}\in\Pi_a$, $y\in U'_{\dot{a}}$, and $z\in U'_{\dot{b}}$. \vspace{3mm}

{\lem\label{lem39}{
For $e\in E(n,\Dt)$ let $\rho(e)=w$ be as in Lemma \ref{lem:38}, and set $e=yx_{{\dot{a}}}(-f)wx_{{\dot{b}}}(g)z$
for $f,g\in D$, ${\dot{a}},{\dot{b}}\in\Pi_a$, $y,\in Y_{\dot{a}}$ and $z\in Y_{\dot{b}}$. Then the following holds:\vspace{3mm}

\underline{Case 1} (for $w_{\dot{a}}(1)e$).\vspace{2mm}

If $f=0$ or $w^{-1}(\dot{a})\in\Delta_a^+$, then $\rho(w_{\dot{a}}(1)e)=w_{\dot{a}}(1)w$.

If $f\neq 0$ and $w^{-1}(\dot{a})\not\in\Delta_a^+$, then $\rho(w_{\dot{a}}(1)e)=h_{\dot{a}}(1)h_{\dot{a}}(f)^{-1}w$.\vspace{3mm}

\underline{Case 2} (for $ew_{\dot{b}}(-1)$)\vspace{2mm}

If $g=0$ or $w(\dot{b})\in\Delta_a^+$, then $\rho(ew_{\dot{b}}(-1))=ww_{\dot{b}}(-1)$.

If $g\neq 0$ and $w(\dot{b})\not\in\Delta_a^+$, then $\rho(ew_{\dot{b}}(-1))=wh_{\dot{b}}(g)h_{\dot{b}}(1)^{-1}$.}}\vspace{3mm}

We now put $X=\{ (e,\tilde{w})\in E(n,\Dt)\times \tilde{N}\ |\ \rho(e)=\psi(\tilde{w})\}$ and define permutations $\lambda(h)$, $\mu(u)$, $\nu_{\dot{a}}$ (resp. $\lambda(h)^*$, $\mu(u)^*$, $\nu^*_{\dot{a}}$) for $\tilde{h}\in\tilde{H}$, $u\in U$ and $\dot{a}\in\Pi_a$ as follows (see Section 2):
\begin{align*}
\lambda(h)(e,\tilde{w})&=(\psi(h)e,h\tilde{w}),\\
(e,\tilde{w})\lambda(h)^*&=(e\psi(h),\tilde{w}h),\\
\mu(u)(e,\tilde{w})&=(ue,\tilde{w}),\\
(e,\tilde{w})\mu(u)^*&=(eu,\tilde{w}),\\
\nu_{\dot{a}}(e,\tilde{w})&=
\begin{cases}
	(w_{\dot{a}}(1)e,\tilde{w}_{\dot{a}}\tilde{w})&\text{if }\rho(w_{\dot{a}}(1)e)=w_{\dot{a}}(1)w,\\
	(w_{\dot{a}}(1)e,\tilde{h}_{\dot{a}}(f)^{-1}\tilde{w})&\text{if }\rho(w_{\dot{a}}(1)e)=h_{\dot{a}}(1)h_{\dot{a}}(f)^{-1}w,
\end{cases}\\
(e,\tilde{w})\nu^*_{\dot{b}}&=
\begin{cases}
	(ew_{\dot{b}}(-1),\tilde{w}\tilde{w}_{\dot{b}}^{-1})&\text{if }\rho(ew_{\dot{b}}(-1))=ww_{\dot{b}}(-1),\\
	(ew_{\dot{b}}(-1),\tilde{w}\tilde{h}_{\dot{b}}(g))&\text{if }\rho(ew_{\dot{b}}(-1))=wh_{\dot{b}}(g)h_{\dot{b}}(1)^{-1},
\end{cases}
\end{align*}
where
\begin{align*}
\tilde{w}_{\dot{a}}&=
\begin{cases}
	\tilde{w}_{\alpha}(1)&\text{if $\dot{a}=(\alpha,0)\in\Pi$},\\
	\tilde{h}_{-\theta}(-t^{-1})\tilde{w}_{-\theta}(1)&\text{if $\dot{a}=\dot{\alpha_0}$},
\end{cases}\\
\tilde{h}_{\dot{a}}(f)&=
\begin{cases}
	\tilde{h}_{\alpha}(f)&\text{if $\dot{a}=(\alpha,0)\in\Pi$},\\
	\tilde{h}_{-\theta}(-f^{-1}t^{-1})\tilde{h}_{-\theta}(-t^{-1})^{-1}&\text{if $\dot{a}=\dot{\alpha_0}$}.
\end{cases}
\end{align*}\vspace{3mm}

\noindent From \cite[Lemma 4.6]{r2},  Lemma \ref{lem:26} and Lemma \ref{lem:210}, we get the following relations.\vspace{3mm}

{\lem\label{lem:310}{
Let $v\in\Dtt$ and $\tilde{w}\in\tilde{N}$. If $\psi(\tilde{w})=P_{\sigma}\mathrm{diag}(u_1,\dots,u_n)$ with some $u_i\in\Dtt$ and the permutation matrix $P_{\sigma}$, then we obtain that
\begin{align*}
&(1)\quad \tilde{w}\tilde{h}_{ij}(v)\tilde{w}^{-1}=\tilde{h}_{\sigma(ij)}(u_ivu_j^{-1})\tilde{h}_{ij}(u_iu_j^{-1}),\\
&(2)\quad \tilde{w}\tilde{w}_{ij}(1)\tilde{w}^{-1}=\tilde{h}_{\sigma(ij)}(u_iu_j^{-1})\tilde{w}_{\sigma(ij)}(1).
\end{align*}
}}

Let $G$ (resp. $G^*$) be the group of automorphisms of $X$ generated by
$\lambda(h)$, $\mu(u)$, $\nu_{\dot{a}}$ (resp. $\lambda(h)^*$, $\mu(u)^*$, $\nu^*_{\dot{a}}$) for all $h\in\tilde{H}$, $u\in U$ and $\dot{a}\in\Pi_a$. We obtain the next lemma as in Section 2.\vspace{3mm}

{\lem[cf. Lemma \ref{lem:214}]\label{lem:311}{
For all $(e,\tilde{w})\in X$, $g\in G$ and $g^*\in G^*$ we have
\[(g(e,\tilde{w}))g^*=g((e,\tilde{w})g^*).\]
}}

{\thm\label{thm:312}{ The group $G$ and $G^*$ operates in a simply transitive manner on $X$, and there exists an epimorphism $\Psi:G\to E(n,\Dt)$ satisfying the following exact sequence which is a central extension:
\[1\to L\to G\xrightarrow{\Psi} E(n,\Dt)\to 1.\]}}
\begin{proof}[PROOF OF THEOREM \ref{thm:32}.]
By Theorem \ref{thm:312} and \cite[Theorem 3]{rs}, we get the following diagram.
\[\begin{split}
\xymatrix{
&1 \ar[r] &L \ar[r]\ar[d]^-{\zeta} &G \ar[r]^-{\Psi}  &E(n,D_{\tau})\ar@{=}[d] \ar[r] &1\\
&1 \ar[r]&K_2(n,\Dt) \ar[r] &St(n,D_{\tau}) \ar[r]\ar[u]^-{p}  &E(n,D_{\tau}) \ar[r] &1
}
\end{split}\]
Then, we see in a way similar to \cite[Proposition 4.9]{r2} that there exists a unique surjective group homomorphism $p: St(n,\Dt)\to G$, so that we find that $p(K_2(n,\Dt))\subset L$ and $p(\hat{c}(u,v))=c(u,v)$ for all $u,v\in\Dt^{\times}$. On the other hand, we know $\zeta$ sends $c(u,v)$ to $\hat{c}(u,v)$ for $u,v\in D_{\tau}$. If we consider the generators of $L$ and $K_2(n,\Dt)$, then we obtain $L\cong K_2(n,\Dt)$.
\end{proof}

{\cor{ We have the following exact sequence:
\[1\longrightarrow K_2(n,\Dt)\longrightarrow Q\stackrel{\varphi}{\longrightarrow} \Dtt\longrightarrow K_1(n,\Dt)\longrightarrow 1.\]}}



\begin{thebibliography}{99}


\bibitem{aabgp}
\textsc{B. Allison, S. Azam, S. Berman, Y. Gao, A. Pianzola},
Extended affine Lie algebras and their root systems,
Mem. Amer. Math. Soc. 126 (1997), no. 603.

\bibitem{b}
\textsc{N. Bourbaki},
\textit{Groupes et alg$\grave{e}$bres de Lie},
Chap. IV-VI, Harmann, Paris (1968).

\bibitem{mj}
\textsc{M. J. Dunwoody},
$K_2$ of a Euclidean ring,
J. Pure Appl. Algebra 7 (1976), no. 1, 53–58.

\bibitem{im}
\textsc{N. Iwahori, H. Matsumoto},
On some Bruhat decomposition and the structure of the Hecke rings of $\mathfrak{p}$-adic Chevalley groups,  Inst. Hautes Études Sci. Publ. Math. No. 25 (1965), 5--48.

\bibitem{mt}
\textsc{R. V. Moody, K. L. Teo},
Tits' systems with crystallographic Weyl groups, J. Algebra 21 (1972), 178--190.

\bibitem{rm}
\textsc{R. Marcuson},
Tits' systems in generalized nonadjoint Chevalley groups, J. Algebra 34 (1975), 84--96.

\bibitem{hm}
\textsc{H. Matsumoto},
Sur les sous-groupes arithm$\grave{e}$tiques des groupes semi-simples d$\grave{e}$ploy$\grave{e}$s (French).
Ann. Sci. École Norm. Sup. (4) 2 (1969), 1–62.

\bibitem{m}
\textsc{J. Milnor},
\textit{Introduction to algebraic $K$-theory},
Ann. of Math. Studies 72, Princeton University Press, Princeton (1971).

\bibitem{jm}
\textsc{J. Morita},
Tits' system in Chevalley groups over Laurent polynomial rings,
Tsukuba J. Math. 3 (1979), no. 2, 41–51.

\bibitem{jmur}
\textsc{J. Morita, U. Rehmann},
A Matsumoto-type theorem for Kac-Moody groups,
Tohoku Math. J. (2) 42 (1990), no. 4, 537–560.

\bibitem{ms}
\textsc{J. Morita, H. Sakaguchi},
Groups defined by extended affine Lie algebras with nullity 2,
Tokyo J. Math. 29 (2006), no. 2, 347–383.

\bibitem{pk}
\textsc{D. H. Peterson, V. G. Kac},
Infinite flag varieties and conjugacy theorems,
Proc. Nat. Acad. Sci. U.S.A. 80 (1983), no. 6, i, 1778--1782.

\bibitem{r2}
\textsc{U. Rehmann},
Zentrale Erweiterungen der speziellen linearen Gruppe eines Schiefk\"{o}rpers (German).
J. Reine Angew. Math. 301 (1978), 77–104.

\bibitem{r}
\textsc{U. Rehmann},
Central extensions of $SL_2$ over division rings and some metaplectic theorems,
Applications of algebraic $K$-theory to algebraic geometry and number theory, Part I, II (Boulder, Colo., 1983), 561–607, Contemp. Math. 55, Amer. Math. Soc., Providence, RI (1986).

\bibitem{hs}
\textsc{H. Sakaguchi},
A Matsumoto-type theorem for linear groups over some completed quantum tori,
Tsukuba J. Math. 32 (2008), no. 1, 1–26.

\bibitem{s}
\textsc{R. Steinberg},
\textit{Lectures on Chevalley groups},
Yale Univ. Lecture Notes, New Haven CT (1968).

\bibitem{rs}
\textsc{R. Sugawara},
Universal central extensions of linear groups over rings of non-commutative Laurent polynomials, associated $K_1$-groups and $K_2$-groups.
Tsukuba J. Math. 45 (2021), no. 1, 13–36.

\bibitem{mt}
\textsc{M. Tomie},
Group presentation of the Schur-multiplier derived from a loop group,
Tsukuba J. Math. 31 (2007), no. 2, 355–395.

\bibitem{yy}
\textsc{Y. Yoshii},
Coordinate algebras of extended affine Lie algebras of type $A_1$,
J. Algebra 234 (2000), no. 1, 128–168.


\end{thebibliography}
\end{document}